\documentclass{iopart}
\usepackage{fontenc}
\usepackage{inputenc}
\usepackage{amssymb}
\usepackage{amsfonts}
\usepackage{stmaryrd}
\usepackage{graphicx}
\usepackage{float}
\usepackage{subfig}
\usepackage{subfloat}
\usepackage{bbold}
\usepackage{hyperref}
\usepackage{csquotes}
\usepackage{bm}
\usepackage{xcolor}

\begin{document}

\title[Coupling rare event algorithms with learned committor functions]{Coupling rare event algorithms with data-based learned committor functions using the analogue Markov chain}

\author{Dario Lucente$^1$, Joran Rolland$^2$, Corentin Herbert$^1$, Freddy Bouchet$^1$}
\address{$^1$ Univ Lyon, ENS de Lyon, Univ Claude Bernard, CNRS, Laboratoire de Physique, F-69342 Lyon}
\address{$^2$ Univ. Lille, Centrale Lille, ENSAM, ONERA, Laboratoire de M\'ecanique des Fluides - Kamp\'e de F\'eriet, UMR 9014, France}

\begin{abstract}
Rare events play a crucial role in many physics, chemistry, and biology phenomena, when they change the structure of the system, for instance in the case of multistability, or when they have a huge impact.
Rare event algorithms have been devised to simulate them efficiently, avoiding the computation of long periods of typical fluctuations.
We consider here the family of splitting or cloning algorithms, which are versatile and specifically suited for far-from-equilibrium dynamics.
To be efficient, these algorithms need to use a smart score function during the selection stage.
Committor functions are the optimal score functions. In this work we propose a new approach, based on the analogue Markov chain, for a data-based learning of approximate committor functions.
We demonstrate that such learned committor functions are extremely efficient score functions when used with the Adaptive Multilevel Splitting algorithm.
We illustrate our approach for a gradient dynamics in a three-well potential, and for the Charney--DeVore model, which is a paradigmatic toy model of multistability for atmospheric dynamics.
For these two dynamics, we show that having observed a few transitions is enough to have a very efficient data-based score function for the rare event algorithm.
This new approach is promising for use for complex dynamics: the rare events can be simulated with a minimal prior knowledge and the results are much more precise than those obtained with a user-designed score function.
\end{abstract}


\maketitle 

\section{Introduction}

Rare events are often extremely important, either because they have a huge impact, like for instance climate extremes~\cite{Seneviratne2021weather},
or because they change completely the structure of the system and shape its history over long times,
like for instance the dynamics of metastability~\cite{farkas1927velocity} and multistability phenomena~\cite{eyring1935activated,kramers1940brownian}.
Such rare events are so important in many physics, chemistry, and biology applications that specific tools have been developed to study them.
These theoretical approaches and dedicated numerical algorithms have been designed by the statistical mechanics and applied mathematics community.

In this paper, we are mainly interested in computational approaches for rare events.
A key difficulty in numerical computation is that these events can be so rare that simulating them directly is prohibitively expensive.
Rare event algorithms and simulations~\cite{bucklew2013introduction}, that aim at reducing their computational cost, have been devised since the nineteen fifties~\cite{kahn1951estimation}.
They have been used to address many problems in statistical physics, for instance studying percolation~\cite{adams2008harmonic},
liquids physics~\cite{allen2017computer}, Lyapunov exponents~\cite{tailleur2007probing}, dynamical phase transitions~\cite{perez2019sampling},
first order phase transitions~\cite{rolland2016computing}, just to cite a few examples among many others.
Chemical physics, biochemistry and the study of biomolecules have inspired many new techniques,
see for example~\cite{bolhuis2002transition,noe2009constructing,metzner2009transition,hartmann2014characterization}.
Recent uses in biology models~\cite{donovan2016unbiased} and ecology have also to be noticed.

Recently, rare events have been studied in far-from-equilibrium systems and non-equilibrium steady states,
where one starts from dynamics without detailed balance.
Rare event techniques have then been extended to scientific fields so far unexpected, with complex dynamics.
For instance in studies of multistability in turbulence~\cite{laurie2015computation,bouchet2019rare},
studies of intermittency in turbulence models~\cite{grafke2013instanton,grafke2015efficient,ebener2019instanton},
transitions to turbulence in pipe and Couette flows~\cite{rolland2018extremely,rolland2021collapse,nemoto2018method},
rogue waves~\cite{dematteis2018rogue}, atmospheric dynamics~\cite{bouchet2019rare,simonnet2021multistability},
climate dynamics~\cite{ragone2018computation,webber2019practical,ragone2020computation,plotkin2019maximizing,finkel2021learning,finkel2020path},
astronomy~\cite{woillez2020instantons,abbot2021rare}, among many other examples.

For such non-equilibrium problems, without detailed balance, one can use either computations related to minimum action methods, possibly related to large deviation theory (see for instance~\cite{grafke2019numerical}),
or the vast family of splitting algorithms or cloning algorithms~\cite{kahn1951estimation,del2012feynman,cerou2007adaptive}.
However, for many applications, for instance in turbulence, climate, atmospheric dynamics, or astronomy, any method that relies on an a priori given bare action is not appropriate.
This is the case when the system is deterministic. This can also be the case for stochastic systems: the precise noise statistics may not be accessible or the rare events may not be produced directly by the model noise but rather by internal fluctuations.

Then, for these cases, the only possible choices for rare event algorithms are splitting algorithms.
These algorithms have indeed been empirically shown to work well for some classes of deterministic chaotic dynamical systems~\cite{wouters2016rare,ragone2018computation}.
An alternative route for studying rare events, without rare event algorithms, would be to use methods that require only short off-equilibrium simulations, for instance through resimulating and milestoning~\cite{noe2009constructing,vanden2009markovian} or coarse-graining of a reduced space of collective variables~\cite{finkel2021learning,finkel2020path}.
Such approaches might be very relevant, however only when the system is simple enough or when one knows sufficiently well the system to define \emph{a priori} relevant collective variables.

The main aim of this paper is to develop the methodology of splitting algorithms such that they might actually be used, practically, for genuinely complex dynamics.
The general principle of splitting algorithms is to perform ensemble simulations, select trajectories prone to produce extremes, discard other less interesting ones,
and resimulate from the interesting ones.
The effectiveness of these algorithms strongly relies on the quality of the score function which is used for the selection stage.
For complex dynamics, in cases when the dynamics is simple enough or the phenomenology of the dynamics is sufficiently well understood to devise good score functions,
splitting algorithms are wonderful tools. For instance, they have been used to compute rare event probabilities
which were totally unreachable with direct numerical simulations, for stochastic partial differential equations~\cite{rolland2016computing},
atmospheric turbulent flows~\cite{bouchet2019rare,simonnet2021multistability}, or full complexity climate models~\cite{ragone2018computation}.
However, without a good score function, splitting algorithms might completely fail.
If the score function is not too bad, but not very good,
splitting algorithms actually produce efficiently rare events, but might suffer from the phenomenon of apparent biases  for the estimation of probabilities~\cite{glasserman1998large,brehier2016unbiasedness}.
The aim of this work is to propose a new methodology to solve these problems and to be able to use splitting algorithms in very complex dynamics without \emph{a priori} knowledge or understanding of a simple effective description of the dynamics.

For many splitting algorithms, there exists a mathematical characterization of an optimal score function: a score function which minimizes the algorithm variance for the computation of the rare event probability and will be very efficient in practice.
For instance, for the Adaptive Multilevel Splitting (AMS)~\cite{cerou2007adaptive}, to be used in this article,
the committor function is the optimal score function~\cite{cerou2019asymptotic}
.
The committor function is the probability that a trajectory visits a region $\mathcal{B}$ of the phase space  before another region $\mathcal{A}$, as a function of the initial condition~\cite{onsager1938initial}
. If $\mathcal{B}$ is the set of rare events of interest, the committor function is then a probabilistic measure of the progress towards the rare event.
The committor function is also a central object of transition path theory~\cite{weinan2005transition,weinan2006towards,vanden2010transition,metzner2006illustration,hummer2004transition}.
A key difficulty is that this optimal score function, the committor function, is actually the rare event probability conditioned on the state of the system.
It contains the information one wishes to compute. One has thus no easy access to it.

\begin{figure*}[tp]
\centering
	\includegraphics[width=15cm]{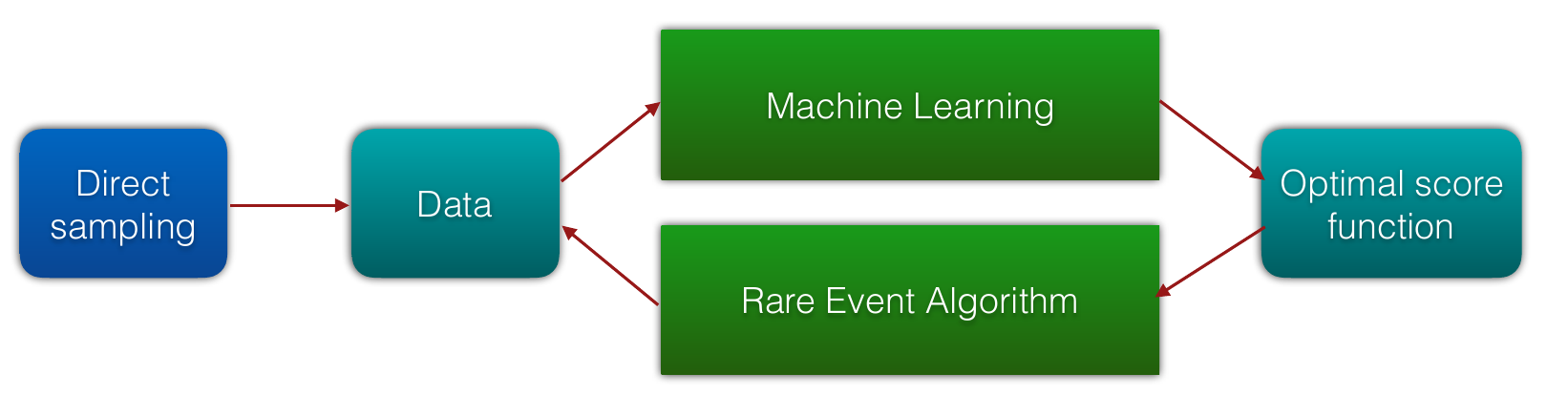} \quad
\caption{Sketch of a feedback control iterative procedure between the rare event algorithm and the machine learning of an approximate optimal score function. The learning of an approximation of the optimal score function makes the algorithm more efficient, and the algorithm provides more data for a better quality of the learning procedure.}
\label{fig:Feedback-Control-Loop}
\end{figure*}

For similar problems, when one would need to know an approximation of a function to efficiently compute the function itself, it is very natural to consider an iterative procedure: a feedback control iterative procedure between the efficient algorithm to produce the data and the learning of the function itself.
The learning of an approximation of the optimal score function makes the algorithm more efficient, and the algorithm provides more data for a better quality of the learning procedure.
This is for instance the idea behind the Wang and Landau algorithm~\cite{wang2001efficient}, in multicanonical methods for equilibrium statistical mechanics, or the idea at the base of adaptive importance sampling~\cite{bugallo2017adaptive}.
This feedback iterative procedure is illustrated in figure~\ref{fig:Feedback-Control-Loop}. We have already implemented such a feedback iterative procedure for the Giardina--Kurchan cloning algorithm, a specific example of a splitting algorithm~\cite{nemoto2016population}. One iteration of the loop was also performed with the AMS algorithm, using a Mondrian forest for learning, for a two-dimensional gradient dynamics~\cite{du2020sequential}. However, the learning step in the first example was extremely simple as the function to be learned was a function over a one-dimensional space. We want to extend this approach to more complex dynamics.

Many interesting methods have been or are currently being devised to learn committor functions: based on direct machine learning~\cite{pozun2012optimizing},
using a characterization of the committor function for diffusions as a solution of a partial differential equation~\cite{khoo2019solving,li2019computing},
computing the committor function from a finite state Markov chain~\cite{schutte1999direct,prinz2011efficient,noe2019markov,tantet2015early},
possibly a Markov state model approximation of the dynamics~\cite{ulam2004problems}.
Recently a very interesting approach has been considered starting from a Galerkin approximation of the dynamics generator, or the Koopman operator.
Finite dimensional approximations of the dynamics generator have been used to identify good reaction coordinates~\cite{froyland2014computational,bittracher2018transition},
or to evaluate eigenfunctions of the operator~\cite{giannakis2015spatiotemporal,giannakis2019data,williams2015data,mardt2018vampnets}, sometimes with climate applications~\cite{giannakis2015spatiotemporal,giannakis2019data}.
Such a direct Galerkin approximation has been used to directly compute committor functions, avoiding the burden of discretizing a high dimensional phase space~\cite{thiede2019galerkin,strahan2021long}.
Several computations of committor functions have been performed with applications in either geophysical fluid dynamics or in climate sciences~\cite{finkel2021learning,miron2021transition,finkel2020path,lucente2020machine,lucente2021committor},
using either direct or involved approaches.

The aim of this paper is to test the coupling of data-based learning of approximate committor functions with rare event algorithms, in the spirit of figure~\ref{fig:Feedback-Control-Loop}. As we are specifically interested in complex dynamics, the learning strategy needs to have the potentiality to scale well in very large dimensions. Moreover, it should be suited for any dynamics, including chaotic deterministic systems or dynamics for which the noise is irrelevant for the process of interest. It also needs to be not too greedy in terms of dataset length.  Among all the possible approaches for learning committor functions, the ones based on approximation of the dynamics generator seem to be best suited~\cite{thiede2019galerkin,strahan2021long}.

In this paper we propose a new method based on an approximation of the generator for the dynamics.  For this purpose, we consider a slightly modified version of the analogue method, first proposed by Lorenz~\cite{lorenz1969three,lorenz1969atmospheric}. The idea behind the analogue method can be summed up by Maxwell's sentence~\cite{garnett1882life} "\textit{From like antecedents follow like consequents}". This approach is nowadays used to build stochastic weather generators~\cite{yiou2014anawege,yiou2019stochastic}. A key remark is that the analogue method defines a Markov chain which is an approximation of the generator of the original dynamics. Then a learned committor function can be computed using classical methods for computing Markov chain committor functions.  This new way to compute committor function, based on the analogue Markov chain, is an alternative path that leads to dynamic-based estimates of the committor function. We show in this paper that this method is actually very simple, robust, and efficient. We show that the learned committor function, based on the analogue Markov chain, is more precise and efficient than the classical $K$-nearest neighbors regression, which computes the committor by averaging the observations of $K$ nearby points.

After having put forward and tested this committor function computation using the analogue Markov chain, we couple it to the Adaptive Multilevel Splitting (AMS)~\cite{cerou2007adaptive}: we directly use the data-based approximate committor function as a score function for the algorithm. We make a precise study that shows that for large enough data sets, the performance of the AMS algorithm is greatly improved. The apparent bias phenomenon is avoided and rare events are computed without \emph{a priori} knowledge of the dynamics.

To summarize the previous discussion, the purpose of this work is twofold. On the one hand, we introduce a data-driven approach which can be used to compute the committor function, and which exploits the dynamical information provided by the observed dynamics. On the other hand, we show how it is possible to use this method to build a learned score function for efficient rare event algorithms. We illustrate our approach for two dynamics. First a stochastic gradient dynamics in a three-well potential, in dimension two. Then we study the Charney--DeVore model, which is a paradigmatic toy model of multistability for atmospheric flows~\cite{charney1979multiple}, with six variables. For these two dynamics, we show that having observed a few transitions is enough to have a very efficient data-based score function for the rare event algorithm.

The paper is organized as follows. In Sec.~\ref{sec:Committor_Function}, we define and discuss the mathematical properties of the committor function, we explain a direct sampling strategy, and define the Brier score which quantifies the quality of an approximate committor function.
Section~\ref{sec:analogue} is devoted to the analogue method and how it can be used to obtain a dynamics-based estimate of the committor function.
Finally, in Sec.~\ref{sec:ams} we introduce the AMS rare event algorithm, we use it with a score function which is the learned analogue Markov chain committor function, and we discuss the improvements given by this approach.

\section{The committor function}\label{sec:Committor_Function}

\subsection{Definition of the committor function for a Markov process}

\begin{figure*}[tp]
\centering
	\includegraphics[scale=.5]{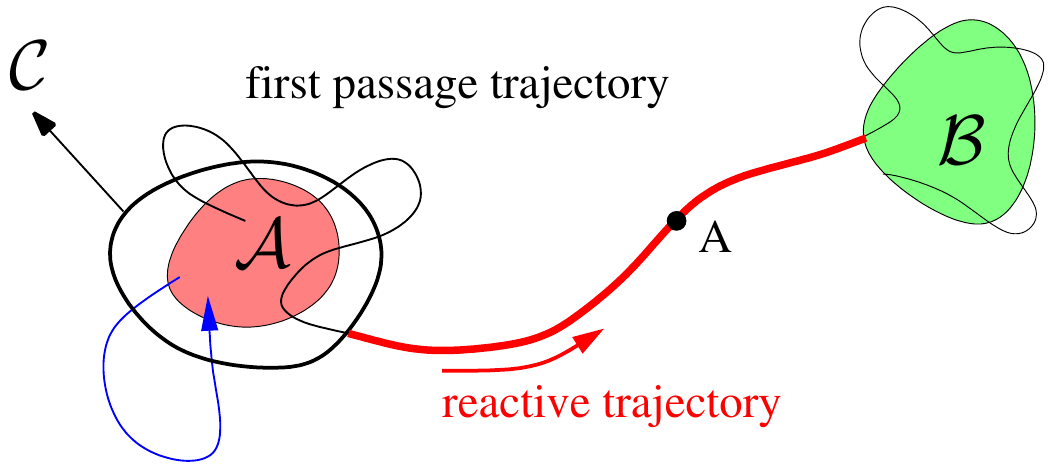} \quad
\caption{Sketch~\cite{rolland2015statistical} of a first passage trajectory from $\mathcal{A}$ to $\mathcal{B}$. The transition path, also called reactive trajectory, is highlighted in red.}
\label{fig:Illustration_reac}
\end{figure*}

For a Markov process, a \emph{committor function}~\cite{weinan2005transition,weinan2006towards,vanden2010transition,metzner2006illustration} is the probability to hit a set $\mathcal{B}$ of the phase space before another set $\mathcal{A}$, conditioned on the knowledge of the initial condition.
In practice, $\mathcal{A}$ and $\mathcal{B}$ may for instance be two regions of phase space corresponding to two metastable states, like two conformations of a protein~\cite{lopes2019analysis} or two configurations of a turbulent flow~\cite{bouchet2019rare,simonnet2021multistability}.
Then, the committor function is the probability of transition from one of these states to the other (see figure~\ref{fig:Illustration_reac}).
Alternatively, region $\mathcal{A}$ may correspond to a typical state around which the system fluctuates, and region $\mathcal{B}$ to an atypical fluctuation of interest because of its impact, usually defined by some observable reaching a given threshold.
In that case, the committor function allows to estimate the probability that the rare event occurs within a given timeframe, or alternatively the return time of the event~\cite{lestang2018computing}.
As an example, we have recently used the committor function within this framework to study the probabilistic predictability of occurrence of El Ni\~no events~\cite{lucente2021committor}.

To give a more precise definition, we consider a discrete time stochastic process on a phase space $\mathcal{X}\subset \mathbb{R}^D$, where $D$ is the space dimension.
A given realization of the process will be noted as ${\{X_n\}}_{1\le n\le N_t}$, with $X_n\in\mathcal{X}$.
The \emph{first hitting time} $T_\mathcal{D}(\mathbf{x})$ of a set $\mathcal{D}\subset\mathcal{X}$ is defined as
\begin{equation}
  T_\mathcal{D}(\mathbf{x})=\textrm{inf}\{n:X_n\in\mathcal{D}|X_0=\mathbf{x}\}.
\end{equation}
The committor function $q(\mathbf{x})$ is the probability that the first hitting time of a set $\mathcal{B}$ be smaller than the first hitting time of set $\mathcal{A}$, as a function of the initial condition, i.e.
\begin{equation}
	q(\mathbf{x})=\mathbb{P}[T_\mathcal{B}(\mathbf{x})<T_\mathcal{A}(\mathbf{x})].
\label{Eq:CommittorDefinition}
\end{equation}
This definition immediately generalizes for continuous time Markov processes.

If the dynamics is a stochastic differential equation, $q(\mathbf{x})$ is the solution of the Dirichlet problem~\cite{weinan2005transition,thiede2019galerkin}:
\begin{equation}
\mathcal{L} q(\mathbf{x})=0 \rm{\ with\ } q(\mathbf{x})=0 \rm{\ if\ } \mathbf{x} \in \mathcal{A} \rm{\ and\ }
 q(\mathbf{x})=1\rm{\ if\ } \mathbf{x} \in \mathcal{B},
 \label{Eq:CommittorDirichletProblem}
\end{equation}
with $\mathcal{L}$ the adjoint of the Fokker-Planck operator:
\begin{equation}
\mathcal{L}=\sum_{i=1}^D{a_i(\mathbf{x})\frac{\partial}{\partial x_i}(\cdot)}+\sum_{i,j=1}^D{\mathbf{b}_{ij}(\mathbf{x})\frac{\partial^2}{\partial x_i\partial x_j}}{(\cdot)},
\end{equation}
where $\mathbf{a}$ is the drift coefficient and $\mathbf{b}$ the diffusion coefficient.
One way to compute a committor function is to solve this partial differential equation. In practice, such a computation is impossible, using standard techniques, as soon as the system has more than a few degrees of freedom. This equation can be used for computing approximate solutions, using machine learning, for systems of dimension $D \sim 10$~\cite{khoo2019solving,li2019computing}.

\subsection{Direct sampling of the committor function}\label{sec:directsampling}

In this section we consider data-based methods for the computation of a committor function. The data consists of sets of trajectories of the stochastic process. The simplest method is to directly use the definition (\ref{Eq:CommittorDefinition}).
In practice, to compute the function at point $\mathbf{x}$, we initialize an ensemble of $N$ trajectories in $X_0=\mathbf{x}$ and evolve them until they reach $\mathcal{A}$ or $\mathcal{B}$.
Let $N_\mathcal{B}$ be the number of trajectories that have reached $\mathcal{B}$ before $\mathcal{A}$.
Then, the value of the committor function at point $\mathbf{x}$ can be estimated as
\begin{equation}
  q(\mathbf{x})=\frac{N_\mathcal{B}}{N}.\label{eq:directsampling}
\end{equation}
Like the Dirichlet problem~\eref{Eq:CommittorDirichletProblem}, this method can only be applied if the equations of motion are known, and it is inapplicable for high dimensional systems, as it requires simulating many trajectories for each point of phase space where we want to compute the committor function. The numerical burden thus increases exponentially with the dimension of the system.

For an ergodic process, the committor function $q(\mathbf{x})$ and the stationary distribution function $\rho(\mathbf{x})$ can be computed from an observed trajectory $\{X_n\}$ from the formulas
\begin{eqnarray}
\rho(\mathbf{x})q(\mathbf{x})=\lim_{N_t\rightarrow\infty}\frac{1}{N_t}\sum^{N_t}_{n=0}\delta\left(X_{n}-\mathbf{x}\right)1_{\{T_{\mathcal{B}}(X_n)\leq T_{\mathcal{A}}(X_n)\}} \  {\rm and} \nonumber \\
\rho(\mathbf{x})=\lim_{N_t\rightarrow\infty}\frac{1}{N_t}\sum_{n=0}^{N_t}\delta\left(X_{n}-\mathbf{x}\right),
\label{Eq:CommittorTrajectory}
\end{eqnarray}
where $\delta$ is a Dirac delta function, and $1_{\{T_{\mathcal{B}}(X_n)\leq T_{\mathcal{A}}(X_n)\}}$ takes value $1$ if the trajectory visits set $\mathcal{B}$ before set $\mathcal{A}$ starting from $X_n$, and $0$ otherwise.
Numerically, $q(\mathbf{x})$ can be computed from~\eref{Eq:CommittorTrajectory} after spatial and temporal discretization of the process (see for instance~\cite{lucente2020machine, lopes2019analysis, lucente2021committor}).
Unlike the previous methods, this approach is applicable even if we do not know the equations of motion. Its numerical cost does not depend on the dimension of phase space, but it only provides estimates of the committor function on points which neighborhood was visited many times by the observed trajectory.

\subsection{Estimating the committor function for any point of the phase space}\label{sec:extendingcommittor}

In Sec.~\ref{sec:directsampling}, we have presented a direct sampling method to estimate the committor function based on data.
However, it provides values only on the set of points that was visited along the trajectory.
This is also true for the other data-based method that we will present in Sec.~\ref{sec:analogue}, the \emph{analogue method}.
For applications, we may need to estimate the value of the committor function for points which were not in the learning dataset.
This may be the case simply for graphical representations of the committor function along a line or on a plane in phase space (e.g.\ Sec.~\ref{sec:example2D}).
Even more importantly, to use the estimated committor function as a score function with the AMS algorithm (Sec.~\ref{sec:ams}), we need to be able to compute it for arbitrary points in phase space.

To do so, we will use a \emph{nearest neighbor method}~\cite{altman1992introduction}.
Let us denote ${\{X_n\}}_{1\le n\le N_t}\in \mathbb{R}^D$ the learning dataset, for which we have an estimate of the committor $\hat{q}(X_n)$.
For any point $\mathbf{y}\in \mathbb{R}^D$, we search the $K$ nearest neighbors (using the Euclidean distance ${d_E(\mathbf{y},\mathbf{x})}^2=\sum_{i=1}^D{(y_{i}-x_{i})}^2$), corresponding to indices $n_j \in \llbracket 1, N_t \rrbracket$ in our dataset, for $1 \leq j \leq K$.
We then perform a weighted average of the corresponding values of the committor:
\begin{equation}
\hat{q}(\mathbf{y})=\frac{\sum_{j=1}^K w_j \hat{q}(X_{n_j}) }{\sum_{j=1}^K w_j}.\label{ker}
\end{equation}
The weights $w_j$ can be chosen uniform: $w_j=1$ (like in Sec.~\ref{sec:committorapplications}) or given by a kernel, such as $w_j=e^{-\frac{{d_E(\mathbf{y},X_{n_j})}^2}{\omega^2}}$, where $\omega > 0$ is a kernel width (like in Sec.~\ref{sec:ams}), depending on the application.

\subsection{Estimation of the quality of an approximate committor function: the Brier score}\label{sec:committorerror}

In this section we address the issue of how to quantify the precision of an estimate of the committor function.
In what follows, the true committor function is denoted by $q$ while $\hat{q}$ stands for our estimate.
As the committor function value $q(\mathbf{x})$ is for any $\mathbf{x}$ the probability of a binary variable, it is natural to look for a score for a forecast of a binary variable.
We also require that this score can be computed directly from observations. The Brier score is a natural candidate.\\

We first consider $Y$ a random variable with binary outcomes, $Y \in \{ 0,1 \}$, and a Bernoulli distribution: $\mathbb{P}[Y=1]=q$ and $\mathbb{P}[Y=0]=1-q$.
In this section $q$ is a single number that does not depend on $\mathbf{x}$.
We look for an estimator that quantifies the precision of an estimation $\hat{q}$ of $q$.

One of the simpler quantities having the required properties was proposed in 1950 by Brier~\cite{brier1950verification}. We consider  $\{Y_n\}_{1\leq n \leq N}$, $N$ independent realizations of the variable $Y$. The Brier score is defined as
\begin{equation}
B_N=\frac{1}{N}\sum_{n=1}^N{(\hat{q}-Y_n)}^2,
\end{equation}
The Brier score is thus a random variable, with values between $0$ and $1$.

The random variable $(\hat{q}-Y_n)^2$ takes value $(1-\hat{q})^2$ with probability $q$ and value $\hat{q}^2$ with probability $(1-q)$.
Then the average value of $B_N$ is
\begin{equation}
\mathbb{E}(B_N) ={(1-\hat{q})}^2q+{\hat{q}}^2(1-q)=q(1-q)+{(\hat{q}-q)}^2.
\label{Eq:BrierAverage}
\end{equation}
The expectation of the Brier score $B_N$ is therefore the sum of two terms. The first one, $q(1-q)$ is related to the stochastic nature of the forecast and is independent of $\hat{q}$ ; it is a fixed lower bound.
Meanwhile, the second term, ${(\hat{q}-q)}^2$, is a quadratic measure of the error made in the estimation of $q$.
The closer the forecast $\hat{q}$ is to the real value $q$, the lower the Brier score is.
While the computation of only the quadratic error requires the knowledge of the truth $q$, the computation of the Brier score does not require the knowledge of $q$.
In the limit $N\to \infty$, we have an ergodic average and $\lim_{N\to\infty}B_N=\mathbb{E}(B_N)$.

We now extend naturally the definition of the Brier score to the case of Markov processes and committor functions, when $q$ is a function that depends of the variable $\mathbf{x}$. We consider a set of events $\{(X_n,Y_n)\}_{1\leq n \leq N}$, where $X_n$ are points in the phase space distributed according to the invariant measure $\rho$ of the Markov process, $\mathbb{E}\left[\delta \left(X_n-\mathbf{x}\right)\right] = \rho(\mathbf{x})$, and $Y_n$ are binary variables which takes the value $1$ with probability $q\left(X_n\right)$ and value $0$ with probability $1-q\left(X_n\right)$. For instance, the couples $(X_n,Y_n)$ can be sampled along one or several trajectories of the Markov chain, where $X_n$ are the states of the Markov chain and $Y_n$ is equal to zero if the first hitting time of $\mathcal{B}$ after $n$ is smaller than the first hitting time of $\mathcal{A}$ after $n$.

We want to estimate the quality of an approximation $\hat{q}$ of the committor function $q$. Then the committor function Brier score is defined as
\begin{equation}
BT_{N}=\frac{1}{N}\sum_{n=1}^{N} {\left[\hat{q}\left(X_n\right)-Y_n\right]}^2,
\label{Eq:BrierFiniteDataSet}
\end{equation}
Extending directly the previous computations, and assuming ergodicity, we have
\begin{equation}
\mathbb{E}(BT_N) = \lim_{N\to\infty}BT_N = \left\| q-\hat{q} \right\|^2_\rho+ \left\| \sqrt{q(1-q)} \right\|^2_\rho ,
\label{Eq:TotalBrier}
\end{equation}
where $\| f \|^2_{\rho}=\int_{\mathcal{D}} f^2(\mathbf{x})\rho(\mathbf{x}) \,{\rm d}\mathbf{x}$ is the $L^2$ norm weighted according to the invariant measure.
Then the committor Brier score is $\left\| q-\hat{q} \right\|^2_\rho$, the weighted $L^2$ norm of the difference $q-\hat{q}$, up to the constant term $\left\| \sqrt{q(1-q)} \right\|^2_\rho$. While the weighted $L^2$ norm cannot be computed without the knowledge of $q$ and $\rho$, the Brier score can be directly computed from the data by the ergodic average (\ref{Eq:BrierFiniteDataSet}).

\section{The analogue Markov chain}\label{sec:analogue}

In this section we introduce the analogue method in one of its current versions~\cite{yiou2014anawege,yiou2019stochastic,lguensat2017analog,platzer2021using,baldovin2018role}.
It provides a way to build effective dynamics from the data that can be reused to generate new trajectories of the system under consideration at a lower computational cost.
Although more precise definitions will be given throughout the section, we think that briefly illustrating the analogue method in its original form proposed by Lorenz~\cite{lorenz1969three,lorenz1969atmospheric} in 1969 is both conceptually and historically instructive.
Furthermore, this can be seen as a particular case of the method we will present in which only $K=1$ analogue is considered.

In a nutshell, the idea is the following.
Suppose we have access to a time series of observations that we will denote by ${\{X_n\}}_{1\le n \le N_t}$, at times $t_n=n\delta t$ where $\delta t$ is the sampling time step. Starting from a state $\mathbf{x}$ at time $t$, we want to predict a possible dynamical evolution after a duration $\Delta t=l\delta t$.
We search among the available data $\{X_n\}_{1\le n \le N_t}$ the closest to $\mathbf{x}$, i.e.\ an analogue, which will be denoted by $X_{n_\star}$:
\begin{equation}
X_{n_\star} = \mathop{\textrm{argmin}}_{\{X_n\}}\{d(\mathbf{x},X_n)\},
\end{equation}
where $d(\cdot,\cdot)$ is a distance.
After identifying the best analogue $X_{n_\star}$, the prediction of $\mathbf{x}(t+\Delta t)$, denoted $\tilde{\mathbf{x}}(t +\Delta t)$, will be
\begin{equation}
\tilde{\mathbf{x}}(t + \Delta t) = X_{n_\star+l}.
\end{equation}
This method was intended by Lorenz as a deterministic prediction. In the following we are rather interested by stochastic predictions, either because the actual dynamics itself is stochastic, or because we understand the analogue method as an approximate effective description of a chaotic dynamics. For stochastic prediction, we will use $K$ analogues rather than a single one.

\subsection{Definition of the analogue Markov chain}\label{sec:analoguedef}

Let ${\{X(t)\}}_{0\le t \le +\infty}$ be a dynamical process that takes values in the phase space $\mathcal{X}\subset \mathbb{R}^D$.
The nature of the process, i.e.\ whether it is deterministic or stochastic, Markovian or not, is irrelevant to the discussion.
Suppose that a realization of this process is observed at regular time intervals $\delta t$ during a total time $\mathcal{T}=N_t\delta t$ and let ${\{X_n\}}_{1\le n\le N_t}$ denote this sampled trajectory made up of $N_t$ points.
Each point $X_n$ is in $\mathbb{R}^D$.

\begin{figure}[htp]
\centering
\includegraphics[width=\linewidth]{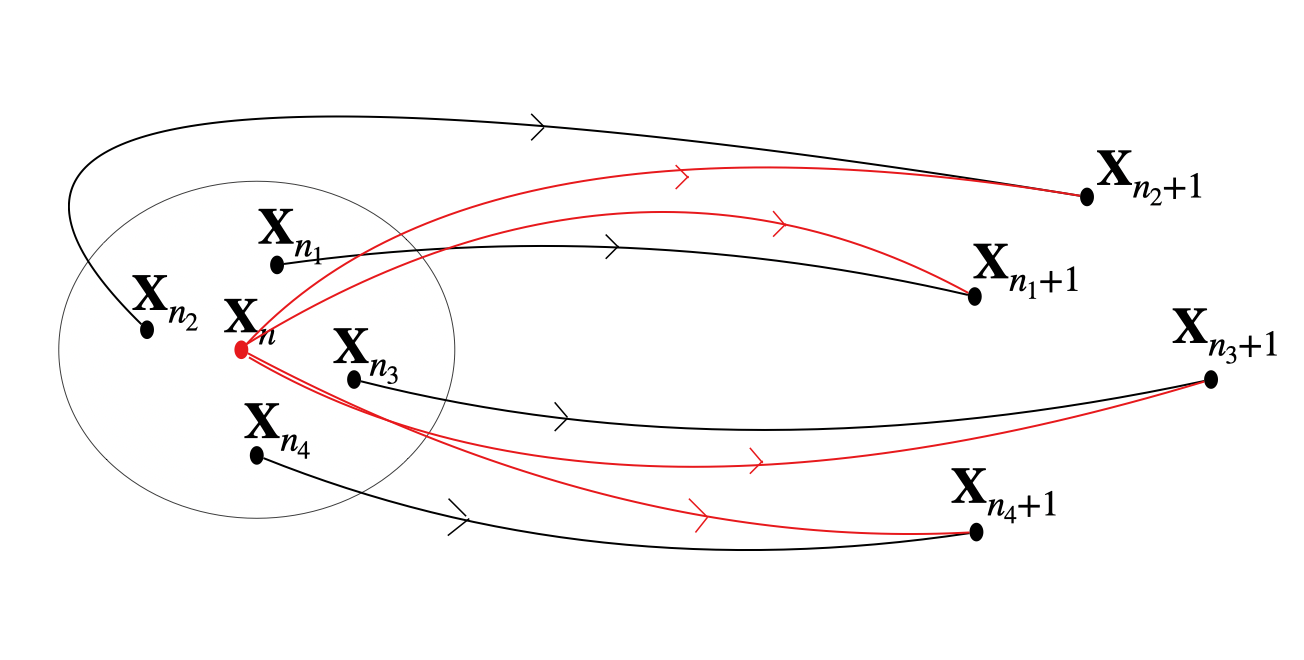}
\caption{Schematic of the analogue method. On the left-hand side of the figure a point $X_n$ surrounded by its analogues $\{\hat{X}_n^{(k)}=X_{n_k}\}_{1 \le k \le K}$ is shown (here $K=4$). On the right-hand side the observed images $\{X_{n_k+1}\}_{1 \le k \le K}$ of the analogues one time step forward are shown. The transitions observed in the data are represented by black lines which link the analogues with their corresponding images. Red lines are associated to the possible transitions from the state $X_n$ of the analogue Markov chain.}
\label{fig:AnaloguesImage}
\end{figure}

We will build a Markov chain that is a data-based approximation of the initial process, based on a generalization of the Lorenz analogue method. We now define possible transitions starting from an observed state $X_n$. Rather than considering just a single nearest neighbor of $X_n$ in the observed data, we will use the $K$ nearest neighbors, where $K$ is a positive number. Those $K$ nearest neighbors are denoted ${\{\hat{X}_n^{(k)}=X_{n_k}\}}_{1 \le k \le K}$, where $n_k \in \llbracket 1, N_t \rrbracket$ is the index of the $k-$th analogue. After identifying analogues $\left\{\hat{X}_n^{(k)}\right\}$, we suppose that we can have a transition between the state $X_n$ and all the possible images of this set of points.
These images will be denoted by $\{X_{n_k+1}\}_{1 \le k \le K}$ and the probability to have a transition between $X_n$ and $X_{n_k+1}$ is set to $\frac{1}{K}$.
An illustration of the analogue Markov chain is shown in figure~\ref{fig:AnaloguesImage}.

With this definition, we see that $K$ is both the number of analogues and the number of possible transitions from any state of the Markov chain. One needs $K$ to be large enough to properly approximate all the possible transitions from a given state. At the same time, the larger $K$, the further the analogue, and the larger the error incurred by using a point further from $X_n$. The optimal value of $K$ will be a tradeoff between these two effects, as a balance between precision and complexity. In practice, $K$ will be chosen empirically, for instance using cross validation. 

\textcolor{black}{It should be noted that, in addition to the number of analogues $ K $, the analogue method depends on a second hyper-parameter, i.e. the lag time $\Delta t$. To lighten the notation we have decided to explain the method for $\Delta t=\delta t$ (where $\delta t$ is the sampling time) but the generalization is straightforward. In analogy to what happens for the choice of $K$, there is no precise protocol for choosing the value of $\Delta t$. In principle, $\Delta t$ should be large enough so that the dynamics of the system on such time scales can be considered Markovian. Hence, $\Delta t$ should be of the same order of magnitude of the correlation time of the system. For the examples in Sec.3, we use $K=150$ and $\Delta t=\delta t$. Although there are no systematic criteria to justify such choices, we stress that for low-dimensional systems, such as those discussed in this paper, the results do not crucially depend on the values of the two hyper-parameters. Therefore, $K$ and $\Delta t$ can be chosen in a fairly wide range. For more complex dynamics, however, further analyses are required to determine the most appropriate values of the hyper-parameters.}

The selection of neighbors is subordinated to the choice of a distance. The best distance most probably depends on the system under investigation.
Distances will be specified on a case-by-case basis.

The analogue Markov chain is a Markov chain on the finite set of $N_t$ observations. In practice, we introduce a matrix with integer entries, $\mathcal{T}\in\mathcal{M}_{N_t K}(\mathbb{N})$.
Each row $n \in \llbracket 1, N_t \rrbracket$ of $\mathcal{T}$ contains the indices of the $K$ nearest neighbors of the point $X_n$, i.e.\ the indices $n_1,\ldots,n_K$ such that $\hat{X}_n^{(k)}=X_{n_k}$ for  $1 \le k \le K$. We stress that $\mathcal{T}$ is not the transition matrix of the Markov chain, to be described latter. $\mathcal{T}$ is rather a matrix of indices of the states.

Since we cannot associate any transition to the end-point $X_{N_t}$, this point will be excluded from the possible candidates for the analogues of each point. To summarize, each entry of $\mathcal{T}$ can take values between 1 and $N_t-1$, i.e. $\mathcal{T}_{nj}\in \llbracket 1, N_t-1\rrbracket$ for all $n,j$ such that $ 1 \leq n \leq N_t $ and $1 \leq j \leq K$.

To generate a synthetic trajectory, we can proceed as follows.
We start with a state $s_0 \in \llbracket 1, N_t \rrbracket$.
Then, we generate a random integer $k$ distributed uniformly in the interval $[1,K]$ and the new state will be $s_1=\mathcal{T}_{s_0k}+1$.
This procedure is iterated to build the entire trajectory.
Through this method we build a Markov chain whose states are ${\{X_n\}}_{1\le n\le N_t}$, i.e.\ the learning dataset.

We now describe the transition matrix $G \in \mathcal{M}_{N_t}(\mathbb{R})$. The elements $G_{nj}$ of $G$ are the probability to observe a transition from the state $n$ to the state $j$. They are given by
\begin{equation}
G_{nj}=\cases{\frac{1}{K} &if $\exists k_{\star} \in \llbracket 1, K \rrbracket: j=\mathcal{T}_{nk_{\star}}+1$,\\
0 &otherwise.\\}
\label{eq:transitionmatrix}
\end{equation}
$G$ is an approximation of the propagator $\mathbb{P}(X_j|X_n)$ of the real dynamics.

Given an observable at time $t$, represented by a column vector $f(t)=f_i(t)$, the observable at time $t+1$ is obtained by applying the operator $G$ to $f(t)$, i.e.
\begin{equation}
f(t+1)=Gf(t).
\end{equation}
Therefore, $G$ plays the same role as the generator of a continuous stochastic process.

Concerning the temporal evolution of probabilities there are two possibilities:
\begin{itemize}
\item consider probabilities as row vectors $\pi$ and let $G$ act to the right, i.e $\pi(t+1)=\pi(t) G$;
\item consider probabilities as column vectors $\pi$ and let them evolve by applying the adjoint operator $G^\dagger$, i.e. $\pi(t+1)=G^\dagger\pi(t)$.
\end{itemize}
In this paper, the second choice has been adopted to emphasize the analogy with continuous stochastic processes.

To initialize a trajectory at a point $\mathbf{x}$ that does not belong to the dataset, we search the $K$ nearest neighbors of $\mathbf{x}$ among the available data and we select as initial condition one of these points with a probability $\frac{1}{K}$. This corresponds to the association of a probability vector $p(\mathbf{x})=p_i(\mathbf{x})$ to the point $\mathbf{x}$ defined as
\begin{equation}
p_{i}(\mathbf{x})=\cases{\frac{1}{K} &if $X_i$ is an analogue of $\mathbf{x}$,\\
0 &otherwise.\\}
\label{eq:initialprobability}
\end{equation}

Note that, for simplicity, in equations~\eref{eq:transitionmatrix} and~\eref{eq:initialprobability} we have assumed that each of the $K$ analogues are chosen with uniform probabilities. We could generalize this choice using analogue dependent weights, for instance computed according to the distances of $X_n$ to its analogues, to account for the varying quality of the various analogues.
The statistical properties of the distance between a state and its different analogues have recently been studied in~\cite{platzer2021using}.

\subsection{Computing the committor function from the analogue Markov chain}\label{sec:analoguecommittor}

Using the analogue Markov chain defined in the previous section, we can compute the committor function $q$ for this Markov chain.
A first approach would be to generate trajectories of this Markov chain, and to directly sample the committor function through a Monte Carlo estimation as described in Sec.~\ref{sec:Committor_Function}.
However, we propose a more efficient computation which consists in solving a linear equation that characterizes the committor function of a Markov chain.
Solving this linear equation is more precise than the direct approach, as we obtain the exact committor function up to numerical accuracy, without sampling errors.
This linear equation will be solved by estimating the leading eigenmodes of a spectral problem, following the algorithm proposed in~\cite{prinz2011efficient}.
Our paper is the first application of this idea to the analogue Markov chain.

We start from the Markov chain transition matrix $G$.
We consider two sets $\mathcal{A} \subset \mathcal{X}$ and $\mathcal{B} \subset \mathcal{X}$, and we will compute the committor function $q$ which is the probability to reach $\mathcal{B}$ before $\mathcal{A}$.
For simplicity, we group together all the states that belong to $\mathcal{A}$ (resp. $\mathcal{B}$) into a single state with index $i_{\mathcal{A}}$ (resp $i_{\mathcal{B}}$).
We then define an auxiliary process where $\mathcal{A}$ and $\mathcal{B}$ are absorbing states: no transition out of these states is allowed.
The corresponding modified transition matrix is $\tilde{G}$, with $\tilde{G}_{i_{\mathcal{A}}i_{\mathcal{A}}}=1$ and for all $j\neq i_{\mathcal{A}}$, $\tilde{G}_{i_{\mathcal{A}}j}=0$, $\tilde{G}_{i_{\mathcal{B}}i_{\mathcal{B}}}=1$ and for all $j\neq i_{\mathcal{B}}$, $\tilde{G}_{i_{\mathcal{B}}j}=0$, while for $i\neq i_\mathcal{A},\,i\neq i_\mathcal{B}$, $\tilde{G}_{ii_\mathcal{A}}=\sum_{k:X_k\in\mathcal{A}}{G_{ik}}$ and $\tilde{G}_{ii_\mathcal{B}}=\sum_{k:X_k\in\mathcal{B}}{G_{ik}}$, and for all other transitions $\tilde{G}_{ij}=G_{ij}$.

For the Markov chain $\tilde{G}$, the committor function is a column vector $q={q_i}$ where $q_i$ is the value of the committor function at the state $i$. $q_i$ is an approximation of the committor function of the initial dynamics at point $X_i$: $q(X_i)$.

For simplicity, we use the same notation for the vector $q$ (associated to the Markov chain) and the function $q$ (associated to the initial dynamics), although they are actually different. In the limit of a large dataset, when the Markov chain fits perfectly the real dynamics, we have asymptotically $q_i \rightarrow q(X_i)$.

From the definition $q_i=\mathbb{P}(T_\mathcal{B}(i)<T_\mathcal{A}(i))$, we have $q_{i_{\mathcal{A}}}=0$ and  $q_{i_{\mathcal{B}}}=1$.
Moreover it is a classical result that $\tilde{G}q=q$~\cite{schutte1999direct,prinz2011efficient,noe2019markov,tantet2015early}.
This is a simple consequence of the estimation of $q$ at two successive steps of the Markov chain.
The affine problem
\begin{equation}
\tilde{G}q=q \qquad {\rm with} \quad q_{i_{\mathcal{A}}}=0  \qquad  {\rm and} \quad  q_{i_{\mathcal{B}}}=1
\label{committor_chain}
\end{equation}
then characterizes the committor function, if we assume that $G$ is ergodic.

Following~\cite{prinz2011efficient}, we note that $1$ is the largest eigenvalue of $\tilde{G}$ (a consequence of the Perron--Frobenius theorem for positive operators that preserve probability).
Moreover $\tilde{G}^{\dagger}$ has two trivial eigenstates with eigenvalue 1, corresponding to situations where the full probability vector is concentrated on state $i_{\mathcal{A}}$ or $i_{\mathcal{B}}$, respectively.
As a consequence, $\tilde{G}$ also has two eigenstates with eigenvalue 1.
If we assume that $G$ is ergodic, then the dimension of the eigenspace of $\tilde{G}$ with eigenvalue 1 is exactly 2.

This gives a simple algorithm to compute $q$.
We first compute $v_1$ and $v_2$ two leading eigenvectors of $\tilde{G}$ with any standard algorithm.
Then $q$ is a linear combination of $v_1$ and $v_2$: $q=\alpha v_1 + \beta v_2$, where $\alpha$ and $\beta$ can be computed from the two conditions $q_{i_{\mathcal{A}}}=0$ and $q_{i_{\mathcal{B}}}=1$.

If the initial dynamics is indeed ergodic, we expect that for large enough dataset the Markov chain $G$ will also be ergodic for most of the realizations.
However, this might not be the case for some realizations.
Such situations could lead to an incorrect computation of $q$ as the solution of~\eref{committor_chain} is then not unique.
In practice we check \textit{a posteriori} (after running the algorithm) whether $q_i\in[0,1]$ for all $i$, which is a necessary condition for $q_i$ to be a probability.
Sometimes, for some realizations of the sampling of the analogue Markov chain, rarely and even more rarely for large datasets, $q$ takes values outside the interval $[0,1]$.
We interpret these cases as a sign of breaking of ergodicity.
We then exclude these rare realizations, with possible ergodicity breaking of the Markov chain, from the results.

\subsection{Applications}\label{sec:committorapplications}

In this section, we estimate the committor function using the analogue method for two different models: Sec.~\ref{sec:example2D} deals with a system of dimension $2$ while Sec.~\ref{sec:charneydevore} concerns a model with $6$ degrees of freedom.
For each system, we compare the estimated committor to the true committor, and we analyze the behavior of the error as the quantity of data upon which the analogue Markov chain relies varies.
Finally, we compare the results of the analogue method with those obtained by the direct method, based on the same amount of data. The committor learned using the analogue Markov chain is denoted by $\hat{q}_A$.

\subsubsection{Model with two degrees of freedom}\label{sec:example2D}

Let us consider a non-trivial $2$-dimensional dynamics~\cite{brehier2016unbiasedness}.
The model is defined by the following stochastic differential equation:
\begin{equation}
\dot{\mathbf{x}}=-\nabla V(\mathbf{x}) +\sqrt{2\epsilon}\mathbf{\Xi}(t),
\label{eq:Dynamics}
\end{equation}
where $\mathbf{x}=(x,y)$, $\mathbf{\Xi}=(\xi_x,\xi_y)$ is a two dimensional gaussian white noise with $\langle \xi_i\rangle =0$, $\langle \xi_i(t)\xi_j(t')\rangle =\delta_{ij}\delta(t-t')$, and the potential $V(\mathbf{x})$ is
\begin{equation}
V(x,y)=0.2x^4+0.2{\left(y-\frac{1}{3}\right)}^4+3e^{-x^2}\left(e^{-{\left(y-\frac{1}{3}\right)}^2}-e^{-{\left(y-\frac{5}{3}\right)}^2}\right)-5e^{-y^2}\left(e^{{(x+1)}^2}+e^{{(x-1)}^2}\right). \label{eq:Potential}
\end{equation}
The stationary distribution of the system is
\begin{equation}
\rho_s(\mathbf{x})=Z^{-1}e^{-\frac{V(\mathbf{x})}{\epsilon}},\label{eq:StationaryDistribution}
\end{equation}
where $Z=\int {\rm d}\mathbf{x}\, e^{-\frac{V(\mathbf{x})}{\epsilon}}$.

\begin{figure}[htbp]
\centering
\subfloat[][\emph{Potential $V(\mathbf{x})$.} \label{fig:Potential}]
	{\includegraphics[width=.45\textwidth]{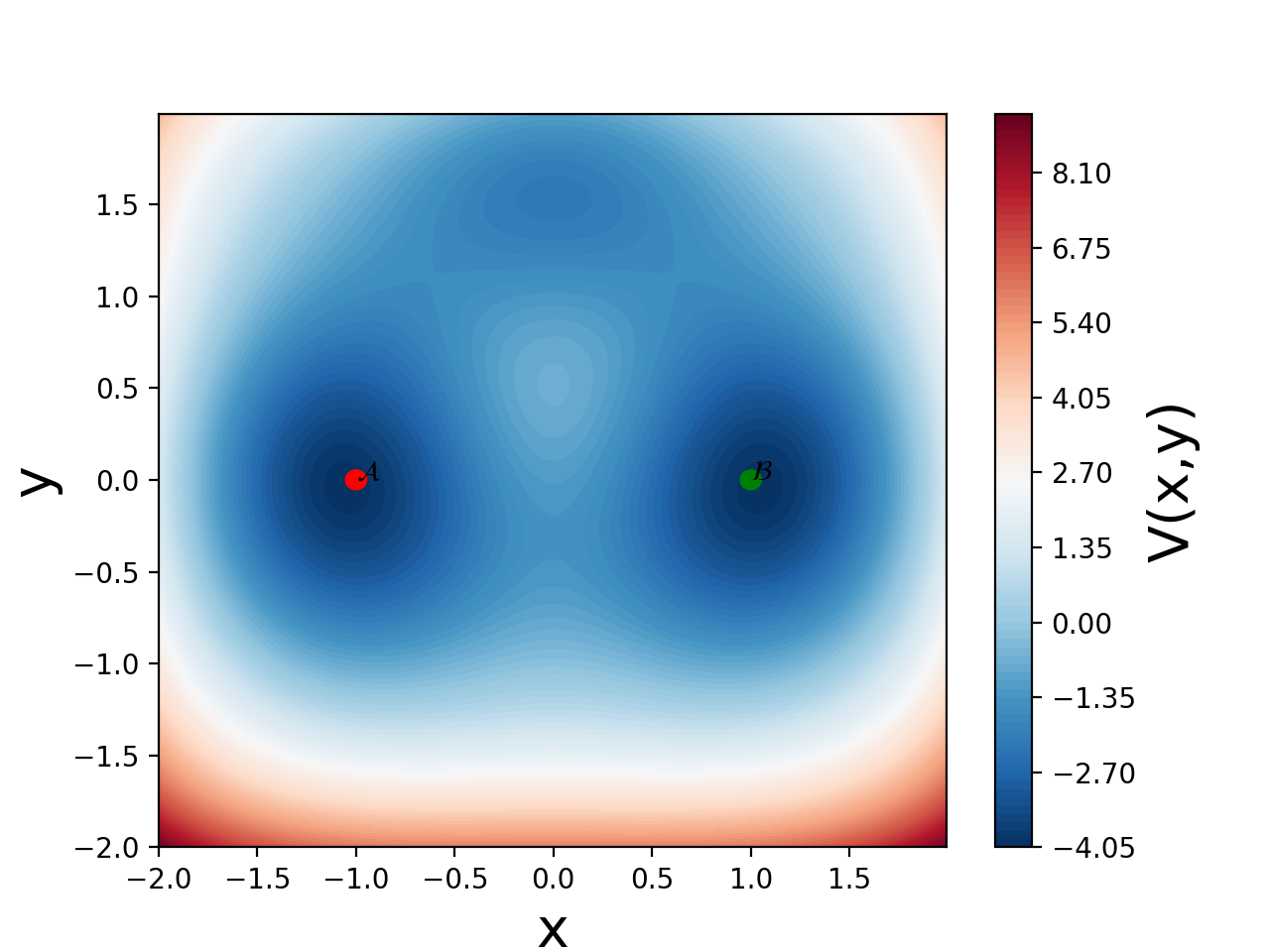} }\quad
\subfloat[][\emph{Stationary distribution $\rho_s(\mathbf{x})$.}\label{fig:TheoreticalDist}]
	{\includegraphics[width=.45\textwidth]{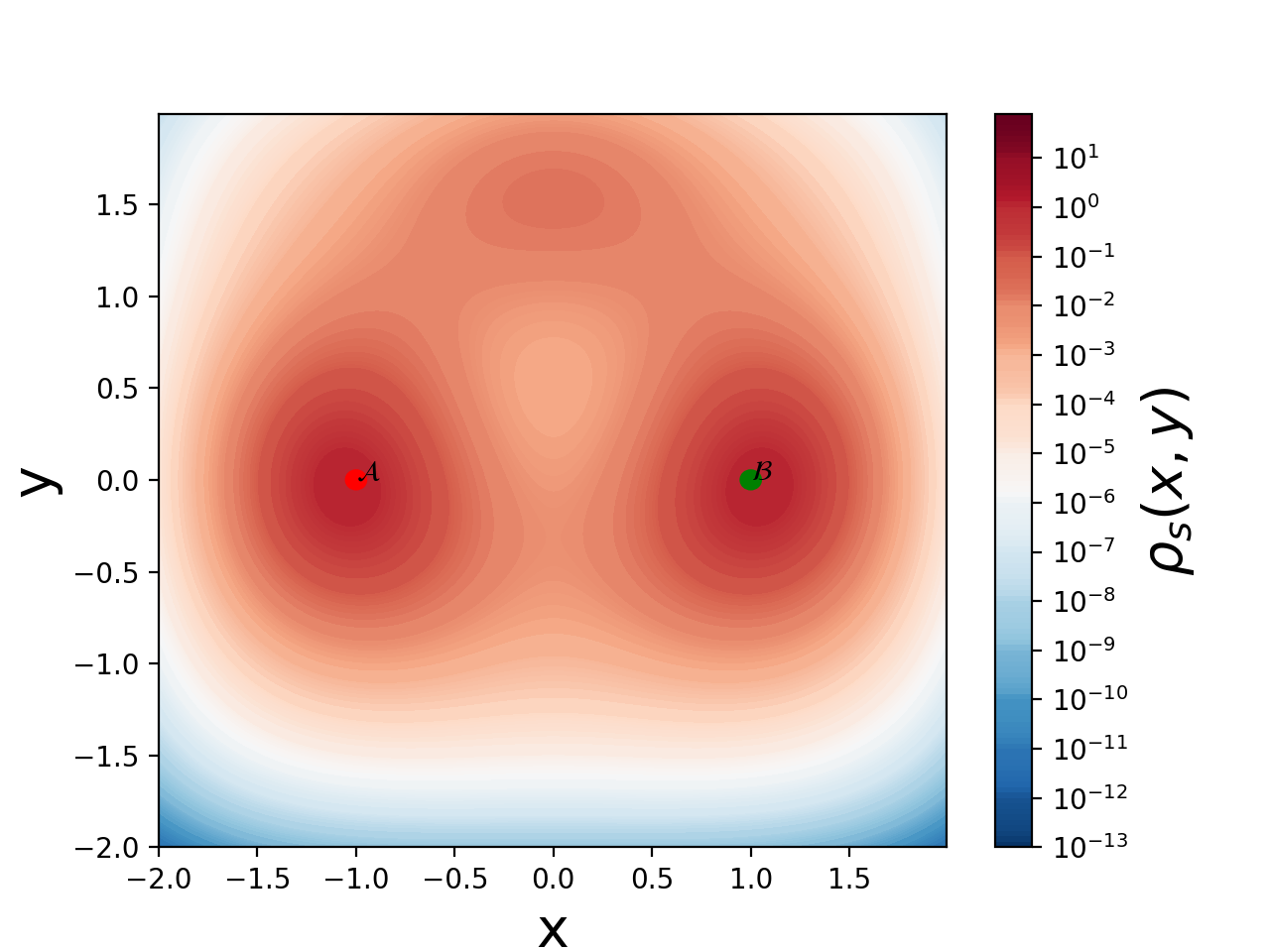}}\quad
\caption{Color maps of the potential $V(\mathbf{x})$ (panel (a), left), defined by~\eref{eq:Potential}, and of the stationary distribution $\rho_s(\mathbf{x})$ (panel (b), right), defined by~\eref{eq:StationaryDistribution}, for $\epsilon=0.5$.}
\label{fig:PotentialAndDistribution}
\end{figure}
Figure~\ref{fig:PotentialAndDistribution} shows both the potential $V(\mathbf{x})$ (\ref{fig:Potential}) and the stationary distribution $\rho_s(\mathbf{x})$ for $\epsilon=0.5$ (\ref{fig:TheoreticalDist}).
As can be seen in figure~\ref{fig:Potential}, $V(\mathbf{x})$ has two global minima close to the points $\mathbf{x}_1=(-1,0)$ and $\mathbf{x}_2=(1,0)$, one local minimum close to the point $\mathbf{x}_m=(0,1.5)$ and a saddle point close to $\mathbf{x}_s=(0,-0.5)$ --- there are also two saddle points separating the global minima from the local minimum, approximately located at $(-0.6, 1.0)$ and $(0.6, 1.0)$.
By comparing the panels~\ref{fig:Potential} and~\ref{fig:TheoreticalDist}, it can be noted that small values of the invariant distribution correspond to large values of the potential and vice versa.
In particular, figure~\ref{fig:TheoreticalDist} shows that $\rho_s(\mathbf{x})$ has global or local maxima at $\mathbf{x}_1$, $\mathbf{x}_2$, and $\mathbf{x}_m$.

Let us consider the two sets $\mathcal{A}=\{\mathbf{x}:d_E(\mathbf{x},\mathbf{x}_1)<0.05\}$ and $\mathcal{B}=\{\mathbf{x}:d_E(\mathbf{x}-\mathbf{x}_2)<0.05\}$.
Note that these sets are defined to include the two maxima of the invariant distribution, where the dynamics spends most of the time.
For $\epsilon=0.5$, the relaxation time $\tau_r$ inside $\mathcal{A}$ or $\mathcal{B}$ is of order $O(1)$, while the average waiting time $T_e$ to observe a transition between these two sets is of order $O(10^2)$.

\begin{figure}[htp]
\centering
\includegraphics[scale=0.6]{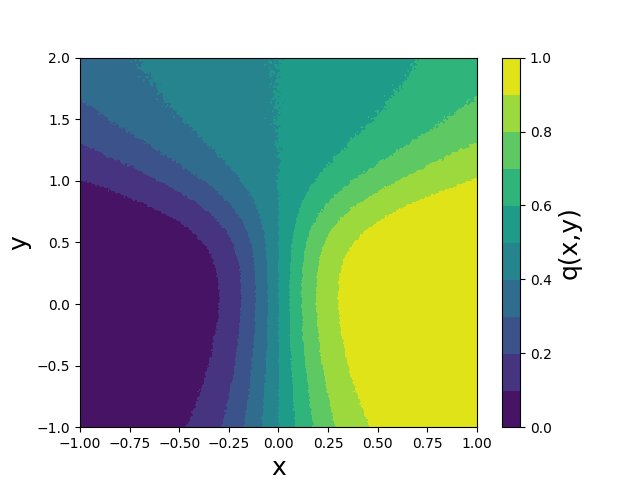}
\caption{Committor function $q(\mathbf{x})$ for the 2D gradient system~\eref{eq:Dynamics}, computed using the real dynamics. The region $x\in [-1,1]$, $y\in [-1,2]$ is divided into $N_{c}=L\times L$ cells $(L=250)$ and, for each cell, $N=10000$ Monte Carlo experiments are performed. }
\label{fig:CommittorReal}
\end{figure}
We will now compute the committor function $q(\mathbf{x})=\mathbb{P}[T_\mathcal{B}(\mathbf{x})<T_\mathcal{A}(\mathbf{x})]$ for this system.
First, we compute a reference committor function in the region $[-1, 1] \times [-1, 2]$ by direct sampling (as explained in Sec.~\ref{sec:directsampling}), using a large amount of data: for each point on a $250\times 250$ grid in this region, we sample 10 000 trajectories until they reach $\mathcal{A}$ or $\mathcal{B}$ and compute the value of the committor at that point using~\eref{eq:directsampling}.
This reference committor function is shown in figure~\ref{fig:CommittorReal}.
One can note that in a region around the set $\mathcal{A}$ the committor function is close to $0$, while in the proximity of $\mathcal{B}$ it is mostly equal to $1$; for $y\simeq -0.5$ and moving along the $x$ direction $q(\mathbf{x})$ changes abruptly through the saddle point $\mathbf{x}_s$.
On the contrary, around the relative minimum point $\mathbf{x}_m$ the committor function is mostly constant, with a value around $0.5$, which corresponds to the probability to reach either of the two minima starting from this point.

\begin{figure*}[tp]
\centering
\subfloat[][\emph{$1$ reactive trajectory.} \label{CommittorGenerator1RT}]
	{\includegraphics[width=.3\textwidth]{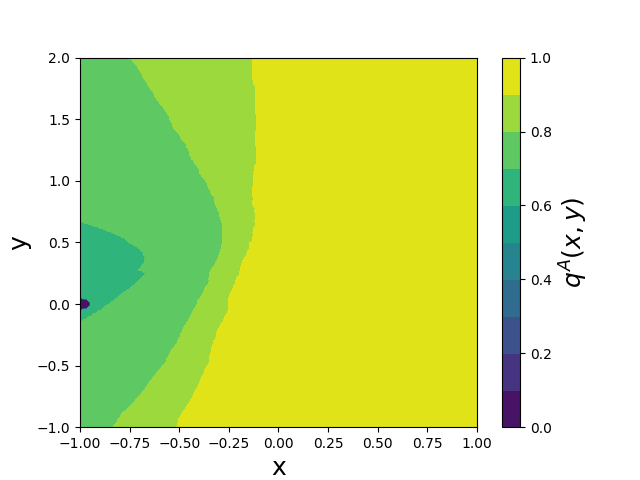} }\quad
\subfloat[][\emph{$2$ reactive trajectories.}\label{CommittorGenerator2RT}]
	{\includegraphics[width=.3\textwidth]{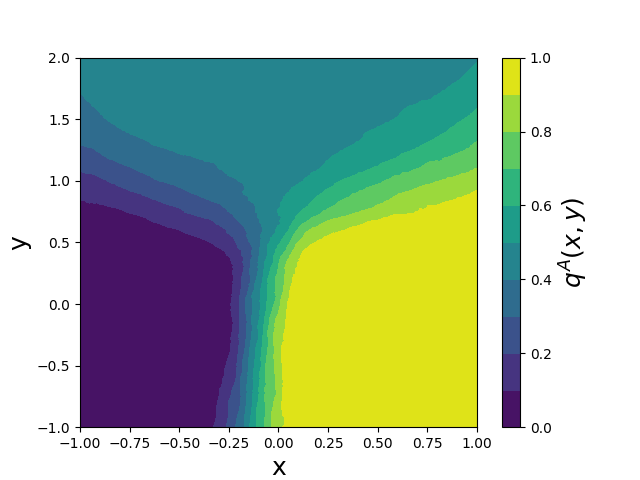}}\quad
\subfloat[][\emph{$20$ reactive trajectories.}\label{CommittorGenerator20RT}]
	{\includegraphics[width=.3\textwidth]{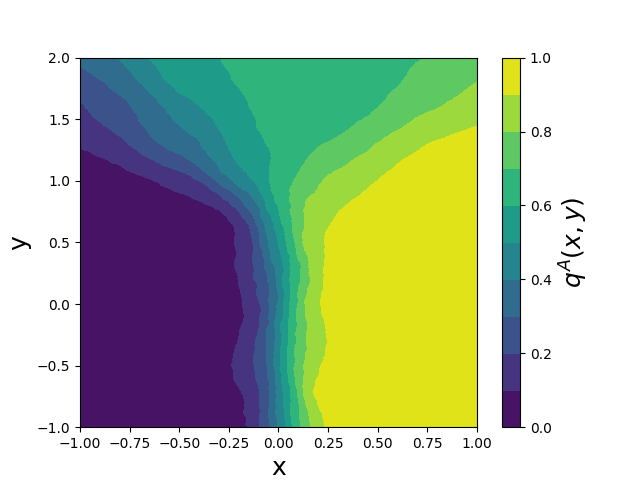}}\\
\subfloat[][\emph{$1$ reactive trajectory.} \label{Committor1RT_NNRegressor}]
	{\includegraphics[width=.3\textwidth]{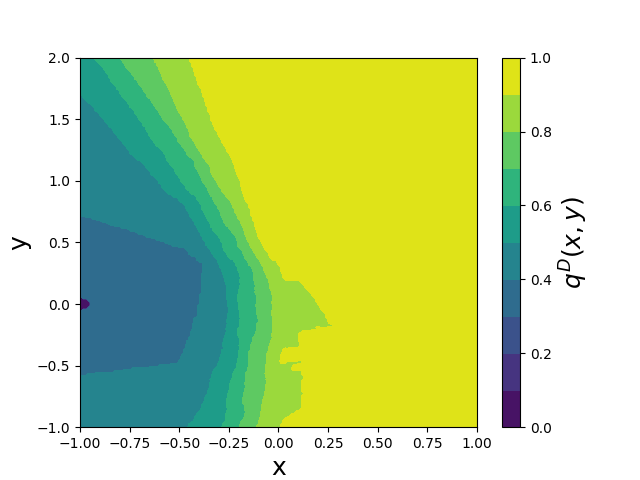} }\quad
\subfloat[][\emph{$2$ reactive trajectories.}\label{Committor2RT_NNRegressor}]
	{\includegraphics[width=.3\textwidth]{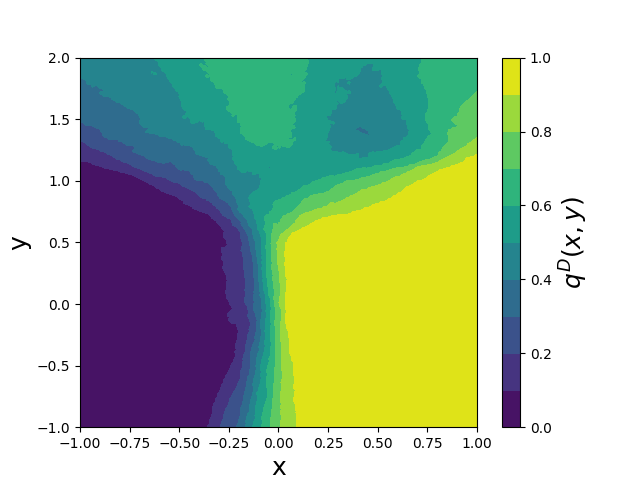}}\quad
\subfloat[][\emph{$20$ reactive trajectories.}\label{Committor20RT_NNRegressor}]
	{\includegraphics[width=.3\textwidth]{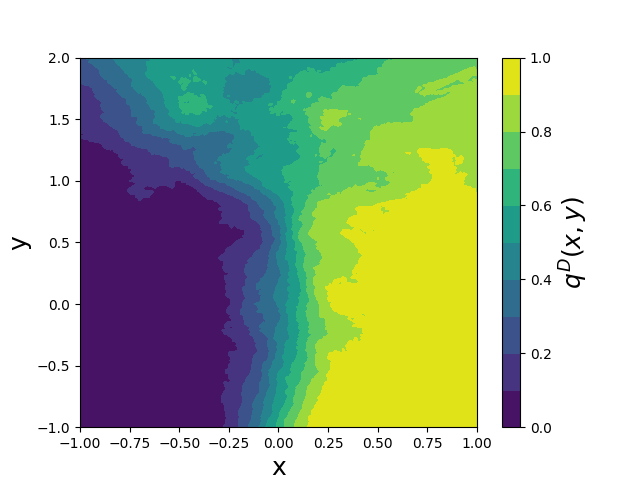}}
        \caption{Estimates $\hat{q}_A(\mathbf{x})$ and $\hat{q}(\mathbf{x})$ of the committor function  of the 2D gradient system~\eref{eq:Dynamics} using the analogue method (top row, panels a--c) and the direct method (bottom row, panels d--f), for learning datasetsets of different length.}
\label{fig:CommittorGeneratorRegressor}
\end{figure*}

Now we estimate the committor function using the method presented in Sec.~\ref{sec:analoguecommittor}.
To do so, we need to generate some learning dataset and to choose a distance and a number of nearest neighbors $K$.
Because we want to measure the quality of our estimator $\hat{q}_A(\mathbf{x})$, by comparing it to the reference committor $q(\mathbf{x})$, as the quantity of available data varies, we generate three trajectories (using the real dynamics) of different length.
Rather than fixing the length of the trajectory, we integrate each trajectory until a fixed number of transitions (1, 2 and 20) between sets $\mathcal{A}$ and $\mathcal{B}$ are observed.
We then construct three analogue Markov chains using each of these trajectories as learning dataset and compute the corresponding committor function.
For these computations, we have used the Euclidian distance and $K=150$ analogues.
The estimate of the committor function for the three choices of learning dataset are shown in Figs.~\ref{CommittorGenerator1RT},\ref{CommittorGenerator2RT},\ref{CommittorGenerator20RT}.
Note that the method presented in Sec.~\ref{sec:analoguecommittor} yields an estimate of the committor function only at the points included in the learning dataset.
To represent the contour levels in Figs.~\ref{CommittorGenerator1RT},~\ref{CommittorGenerator2RT} and~\ref{CommittorGenerator20RT}, we extend our estimate of the committor function to the whole region of interest by using a \emph{$K$-nearest neighbor regression} method, as explained in Sec.~\ref{sec:extendingcommittor}.
To avoid introducing additional parameters, we choose uniform weights $w_j=1$ for all the nearest neighbor and we use the same number of neighbors as for constructing the analogue Markov chain \textbf{$K=150$}.

In addition to the reference committor, we also want to compare the committor estimator based on the analogue method to a direct sampling estimate with the same amount of data.
To do so, we also compute the committor function using~\eref{Eq:CommittorTrajectory} for the same three trajectories as above.
In practice, because the exact same points are never visited twice, this amounts to assigning value 1 to a point in the trajectory if set $\mathcal{B}$ is visited before $\mathcal{A}$ in the rest of the trajectory, and value 0 otherwise.
Again, this provides an estimate of the committor function only at points included in the learning dataset and we extend it to the region of interest with the same $K$-nearest neighbor method as above.
This alternative estimator for the committor function, which we refer to as the \emph{direct method}, is shown in Figs.~\ref{Committor1RT_NNRegressor},~\ref{Committor2RT_NNRegressor}, and~\ref{Committor20RT_NNRegressor}.

Several conclusions can be drawn by comparing qualitatively the committor estimates shown in figure~\ref{fig:CommittorGeneratorRegressor} with the reference committor shown in figure~\ref{fig:CommittorReal}.
First of all, note that a single reactive trajectory does not contain enough information to capture the structure of the committor function (Figs.~\ref{CommittorGenerator1RT},\ref{Committor1RT_NNRegressor}).
The committor estimates start to be qualitatively acceptable when two reactive trajectories are used (Figs.~\ref{CommittorGenerator2RT},\ref{Committor2RT_NNRegressor}).
This is due to the fact that our data set includes the two types of transition paths between $\mathcal{A}$ and $\mathcal{B}$ (the one that passes through the saddle point $\mathbf{x}_s$ and the one that goes through the relative minimum $\mathbf{x}_m$).
A comparison between figure~\ref{CommittorGenerator2RT} and figure~\ref{Committor2RT_NNRegressor} shows that the analogue method gives smoother results than the direct approach.
However, note that both methods have a sharper transition region than that shown in figure~\ref{fig:CommittorReal}.
By increasing the number of reactive trajectories, a wider transition region is obtained (see Figs.~\ref{CommittorGenerator20RT},\ref{Committor20RT_NNRegressor}) and the results appear more similar to the reference committor.
Again, note that the result of figure~\ref{CommittorGenerator20RT} is smoother than that of figure~\ref{Committor20RT_NNRegressor}.

To quantify the error made in approximating $q(\mathbf{x})$ we consider the quantity
\begin{equation}
\| q-\hat{q}\|^2_{\rho_s}=\int{\rm d}\mathbf{x} {[q(\mathbf{x})-\hat{q}(\mathbf{x})]}^2\rho_s(\mathbf{x})\approx \frac{1}{N_p}\sum_{i=1}^{N_{p}}{(q(\mathbf{x}_i)-\hat{q}(\mathbf{x}_i))}^2,
\label{Eq:ErrorEstimation}
\end{equation}
where $q$ is the true committor and $\hat{q}$ the estimate.
The justification of using~\eref{Eq:ErrorEstimation} as an error measurement has already been given in Sec.~\ref{sec:committorerror}: $\| q-\hat{q}\|^2_{\rho_s}$ corresponds to the non-constant term in the Brier score~\eref{Eq:TotalBrier}.
\begin{figure}[tp]
\centering
\includegraphics[width=.45\textwidth]{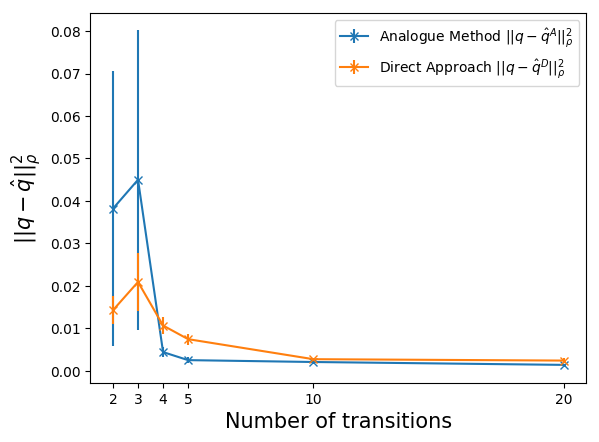}
\caption{Error for the analogue and direct estimators of the committor function of the 2D gradient system~\eref{eq:Dynamics} as function of the number of reactive trajectories in the learning dataset.}
\label{fig:Error_Committor}
\end{figure}
The errors computed from~\eref{Eq:ErrorEstimation} for the two estimators of the committor function (the analogue method and the direct method) are represented in figure~\ref{fig:Error_Committor} as a function of the number of transitions in the dataset.
Each point in figure~\ref{fig:Error_Committor} corresponds to the average error computed over $10$ independent realizations containing the same number of transitions while the error bar corresponds to the standard deviation.
It can be noted that, for small datasets ($2$ or $3$ reactive trajectories), the performances of the two methods are comparable within statistical errors but the direct approach seems to provide more stable results.
A simple interpretation is that when there is not enough data, the analogue Markov chain is not a good enough approximation of the real dynamics to provide any benefit to estimate the committor function.
However, it becomes the case as the amount of data increases, and the analogue method outperforms the direct method as soon as the learning dataset contains at least 4 transitions. When the data contains at least 4 transitions, the error with the analogue method is two to three times smaller than the error with the direct method.

\subsubsection{The Charney--DeVore model}\label{sec:charneydevore}

We now apply the analogue method to compute a committor function for a more complex dynamics, the Charney--DeVore model~\cite{charney1979multiple}.
It is a simple toy model of atmospheric dynamics in the Northern Atlantic region, represented as a 2D channel with differential rotation.
This model was introduced with the aim of proving that the combination of topography and barotropic instabilities can lead to different atmospheric flow regimes.
 It is not intended to be realistic. Actually, the kind of multistability observed in this model is not observed in real atmospheric dynamics. The interest of this model is more methodological, providing a relevant dynamics of intermediate complexity.
The model is obtained by expanding the quasi-geostrophic stream function $\psi(z, y, t)$ ($z$ corresponds to the longitude and $y$ to the latitude) on the basis $\{\phi_{nm}(z,y)\}$ with
\begin{eqnarray}
    \phi_{0m}=\sqrt{2}\cos{\Big(\frac{my}{b}\Big)},\\
    \phi_{nm}=\sqrt{2}\exp{(inz)}\sin{\Big(\frac{my}{b}\Big)},
\end{eqnarray}
and truncating the series to retain only the first six terms.
After the following change of variables~\cite{de1989analysis},
\begin{eqnarray}
  x_1&=\frac{1}{b}\psi_{01}, &x_4=\frac{1}{b}\psi_{02},\\
  x_2&=\frac{1}{\sqrt{2}b}(\psi_{11}+\psi_{-11}), &x_5=\frac{1}{\sqrt{2}b}(\psi_{12}+\psi_{-12}),\\
  x_3&=\frac{i}{\sqrt{2}b}(\psi_{11}-\psi_{-11}), &x_6=\frac{i}{\sqrt{2}b}(\psi_{12}-\psi_{-12}),
\end{eqnarray}
the truncated equations of motion become
\begin{eqnarray}
  \dot{x_1}=\tilde{\gamma}_1x_3-C(x_1-x_1^\star)+\sqrt{2\epsilon}\xi_1\,, \nonumber\\
  \dot{x_2}=-(\alpha_1x_1-\beta_1)x_3-Cx_2-\delta_1x_4x_6+\sqrt{2\epsilon}\xi_2\,,\nonumber\\
  \dot{x_3}=(\alpha_1x_1-\beta_1)x_2-\gamma_1x_1-Cx_3+\delta_1x_4x_5+\sqrt{2\epsilon}\xi_3\,,\nonumber\\
  \dot{x_4}=\tilde{\gamma}_2x_6-C(x_4-x_4^\star)+\eta(x_2x_6-x_3x_5)+\sqrt{2\epsilon}\xi_4\,,\nonumber\\
  \dot{x_5}=-(\alpha_2x_1-\beta_2)x_6-Cx_5-\delta_2x_3x_4+\sqrt{2\epsilon}\xi_5\,,\nonumber\\
  \dot{x_6}=(\alpha_2x_1-\beta_2)x_5-\gamma_2x_4-Cx_6+\delta_2x_2x_4+\sqrt{2\epsilon}\xi_6\,, \label{eq:CharneyDeVoreModel}
\end{eqnarray}
where a Gaussian white noise $\boldsymbol{\xi}(t)$ has been added with an arbitrary amplitude controlled by the parameter $\epsilon$.
All the components of the noise are independent and delta-correlated in time: $\langle\xi_i(t)\xi_j(t')\rangle=\delta_{ij}\delta(t-t')$.
The parameters in (\ref{eq:CharneyDeVoreModel}) are defined as follows
\begin{eqnarray}
    \alpha_m&=\frac{8\sqrt{2}}{\pi}\frac{m^2}{4m^2-1}\frac{b^2+m^2-1}{b^2+m^2}, & \tilde{\gamma}_m=\gamma\frac{4m}{4m^2-1}\frac{\sqrt{2}b}{\pi},\\
    \beta_m&=\frac{\beta b^2}{b^2+m^2}, & \eta=\frac{16\sqrt{2}}{5\pi},\\
    \delta_m&=\frac{64\sqrt{2}}{15\pi}\frac{b^2-m^2+1}{b^2+m^2}, & \gamma_m=\gamma\frac{4m^3}{4m^2-1}\frac{\sqrt{2}b}{\pi(b^2+m^2)}.
\end{eqnarray}
There are $7$ free parameters in this model: $b,\gamma,\beta,C,x_1^\star,x_4^\star$, and the noise amplitude $\epsilon$.
For $\epsilon=0$, the main feature of the system is the coexistence of multiple equilibrium states, in particular the existence of blocked flow and zonal flow regimes.
The number and stability of these equilibrium states depend on the choice of the system parameters~\cite{de1989analysis,crommelin2004mechanism}.
We adopt the same choice made by  Grafke \emph{et al.}~\cite{grafke2017long,grafke2019numerical}, that is $\{b,\gamma,\beta,C,x_1^\star,x_4^\star\}=\{0.5,1,1.25,0.1,4.5,-1.8\}$.
Crommelin \emph{et al.} show that for these parameter values the system has two stable equilibrium points~\cite{crommelin2004mechanism}: one corresponding to a zonal regime and the other to a blocked one.
\begin{figure*}[tp]
\centering
\subfloat[][\emph{Convergence to the zonal regime.} \label{ConvergenceToZonal}]
	{\includegraphics[width=.45\textwidth]{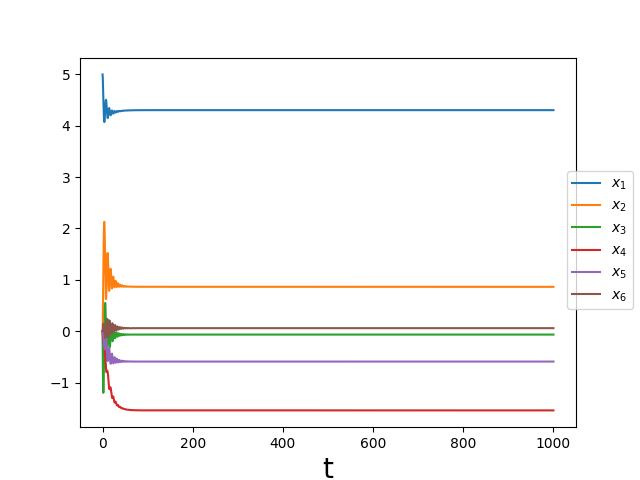} }\quad
\subfloat[][\emph{Convergence to the blocked regime.}\label{ConvergenceToBlocked}]
	{\includegraphics[width=.45\textwidth]{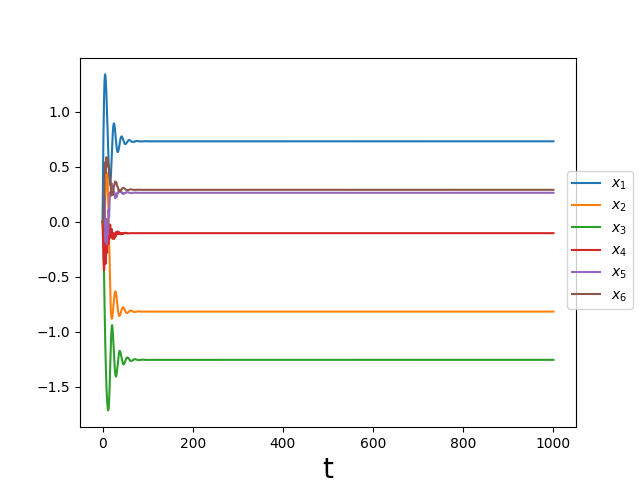}}\\
\subfloat[][\emph{Stream function $\psi(z,y,t)$ in the zonal regime.} \label{StreamZonal}]
	{\includegraphics[width=.45\textwidth]{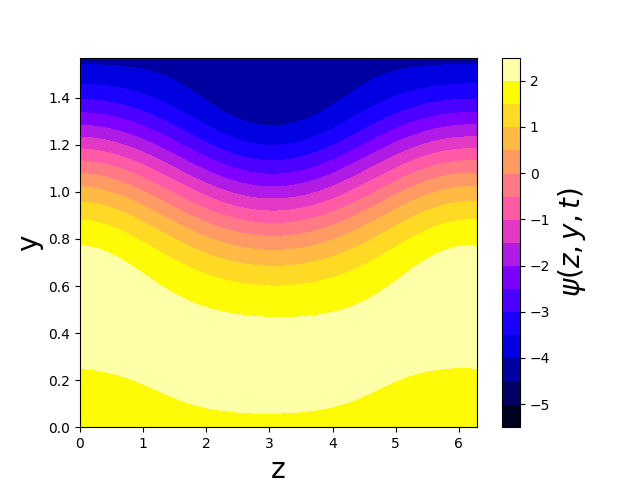} }\quad
\subfloat[][\emph{Stream function $\psi(z,y,t)$ in the blocked regime.}\label{StreamBlocked}]
	{\includegraphics[width=.45\textwidth]{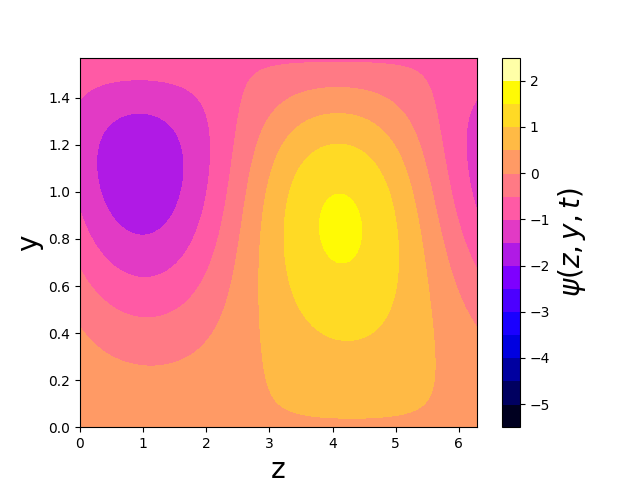}}\\
\caption{Time evolution of the six variables of the Charney--DeVore model for different initial conditions ((a) and (b)) showing the relaxation to two different states and the corresponding stream functions ((c) and (d)).}
\label{fig:ZonalBlocked}
\end{figure*}
Figure~\ref{fig:ZonalBlocked} shows the convergence of the system towards the two equilibrium states as well as the corresponding stream function $\psi$ for the deterministic model ($\epsilon=0$).
The panels~\ref{ConvergenceToZonal} and~\ref{ConvergenceToBlocked} show that, for this choice of parameters, the system exhibits multistability, and that the time it takes to reach the stationary regimes is of order $O(10)$.
The two equilibria correspond to a \emph{zonal} state, with almost horizontal streamlines (figure~\ref{StreamZonal}) and a \emph{blocked} state, with strong cyclonic and anticyclonic structures (figure~\ref{StreamBlocked}).
In the zonal regime the flow is characterized by a strong eastward jet $u_z=\partial_y \psi(z,y,t)$.
Instead, in the blocked state there is no jet, the flow meanders strongly across the domain and it is characterized by the presence of vorticity.

For $\epsilon\neq0$, the system can switch spontaneously from one regime to the other, under the influence of noise.
To study the noise-induced transitions between the zonal and blocked states, we have to define the corresponding regions of the phase space.
Let $\mathbf{x}_{eq}^Z$ and $\mathbf{x}_{eq}^B$ be the equilibrium points corresponding to zonal and blocked flow, respectively.
Given two radii $r_B,r_Z>0$, we define the sets
\begin{eqnarray}
    \mathcal{A}=\{\mathbf{x}\,:\,d_E(\mathbf{x},\mathbf{x}_{eq}^Z)<r_Z\}\,,\nonumber\\
    \mathcal{B}=\{\mathbf{x}\,:\,d_E(\mathbf{x},\mathbf{x}_{eq}^B)<r_B\}\,.
    \label{eq:setscdv}
\end{eqnarray}
In the rest of this section, we consider $r_Z=0.8$, $r_B=0.3$ and $\epsilon=0.02$.
For such parameters, the average time between two transitions is of order $O(10^3)$.

Let us now discuss the committor function $q(\mathbf{x})=\mathbb{P}(T_\mathcal{B}(\mathbf{x})<T_\mathcal{A}(\mathbf{x}))$ of the system.
First of all, it should be noted that a direct computation of $q(\mathbf{x})$ in the whole phase space is not feasible.
Indeed, such a calculation would require discretizing the six-dimensional phase space and to simulate a set of $N$ trajectories for each point of the domain until they reach either $\mathcal{A}$ or $\mathcal{B}$.
If $100$ points along each direction were to be taken, then $N\times10^{12}$ trajectories would have to be simulated.
Considering a time of one millisecond to simulate $N$ trajectories, the computation of $q(\mathbf{x})$ would still take $T_q=10^9\,\,\textrm{s}\approx 11574\,\,\textrm{days}$.
Therefore, the reference committor $q(\mathbf{x})$ is computed on a limited number of points $N_p$ distributed according to the invariant measure.
Since the invariant distribution of the system is not known, the points $N_p$ are sampled at regular time intervals over a very long trajectory.
To be more specific, we consider a trajectory $10^7$ time units long and we sample the $N_p$ points at intervals $\delta t=10^3$ time units.
In this way, we ensure the statistical independence of the points. Furthermore, their distribution will coincide with the invariant distribution of the system in the limit $N_p\to + \infty$ by construction.
Then, the committor function on those points can be computed by running $N$ Monte Carlo experiments for each of them.

After computing $q(\mathbf{x})$ along a trajectory in the six-dimensional space, it is natural to ask how to represent such a function in a low-dimensional space.
We will show below the empirical distribution of the value of the committor conditioned on the coordinate $x_1$, defined as:
\begin{equation}
  \zeta(q|x_1) = \frac{\int{\rm d}\mathbf{y}\rho_s(\mathbf{y})\delta(q(\mathbf{y})-q)\delta(y_1-x_1)}{\int{\rm d}\mathbf{y}\rho_s(\mathbf{y})\delta(y_1-x_1)}.
\end{equation}
We chose to condition on the coordinate $x_1$ because the separation of the two attractors is larger in this direction than in the others. \textcolor{black}{We want to emphasize that the distribution $\zeta(q|x_1)$ is introduced for illustrative purposes only. We stress that we actually compute the committor function $q(\mathbf{x})$ on the whole phase space, and we use all these values to evaluate the performance of the analogue method. Similarly, the score function that we use in the rare event algorithm is the committor function $q(\mathbf{x})$.}
\begin{figure}[htp]
\centering
\includegraphics[width=0.5\textwidth]{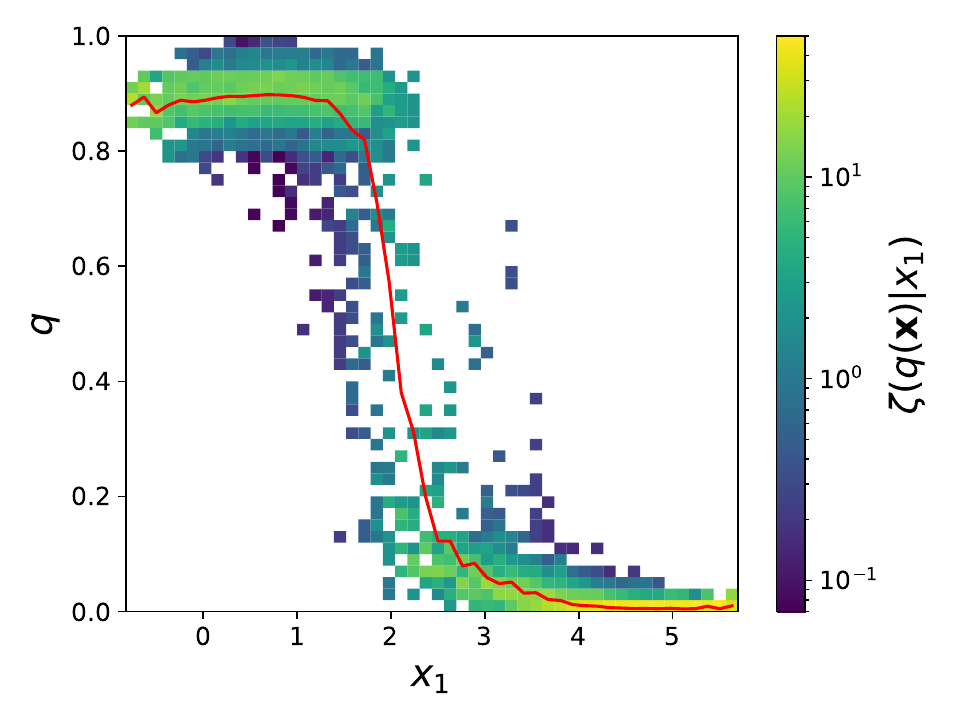}
\caption{Color map of the conditional distribution $\zeta(q|x_1)$ for the reference committor of the Charney--DeVore model. The conditional average $\langle q\rangle_{x_1}$ as a function of $x_1$ is shown as a red line.
The reference committor is computed using the real dynamics on $N_p=10000$ points of the phase space distributed according to the stationary distribution $\rho_s$.
For each point, $N=100$ Monte Carlo experiments are performed.}
\label{fig:CommittorReferenceN100}
\end{figure}
The distribution $\zeta(q|x_1)$ for the reference committor is represented on figure~\ref{fig:CommittorReferenceN100},
along with the conditional average $\langle q\rangle_{x_1}=\int q\zeta(q|x_1)\,{\rm d}q $.
We can first consider the conditional average of $q$ as a function of $x_1$.
It is close to $1$ for $x_1\lesssim 1.5$ as $x$ is close to $\mathbf{x}_{eq}^B$.
The conditional average of the committor decreases like a sigmoid for $1.5\lesssim x_1\lesssim 2.5$.
It is finally close to $0$ for $x_1\gtrsim 2.5$.
Note however that the conditional average misses a lot of information in the range $1.5\lesssim x_1\lesssim 2.5$ that the conditional distribution give us.
The conditional distribution is extremely dispersed in this range.
Actually, in this range the conditional distribution is bimodal and considering the committor as a function of $x_1$ does not give us all the information on this function.
The fact that the committor function exhibits some spread around 0 and 1 close to sets $\mathcal{A}$ and $\mathcal{B}$ can be explained by observing that many of the points are in fact located outside of the hyperballs defining these sets~\eref{eq:setscdv}, although they lie in the basins of attraction of the zonal and blocked states.

We now estimate the committor function for the Charney--DeVore model using the analogue method.
As in Sec.~\ref{sec:example2D}, we use several datasets of different size to build the analogue Markov chain used to estimate the committor.
The size of these datasets is measured by the number $n=2$, 5, 10, 15 of transitions between $\mathcal{A}$ and $\mathcal{B}$.
As previously, we select $K=150$ analogues using the Euclidean distance.
\begin{figure*}[tp]
\centering
\subfloat[][\emph{$2$ reactive trajectory.} \label{CommittorAnalogue_2}]
	{\includegraphics[width=.48\textwidth]{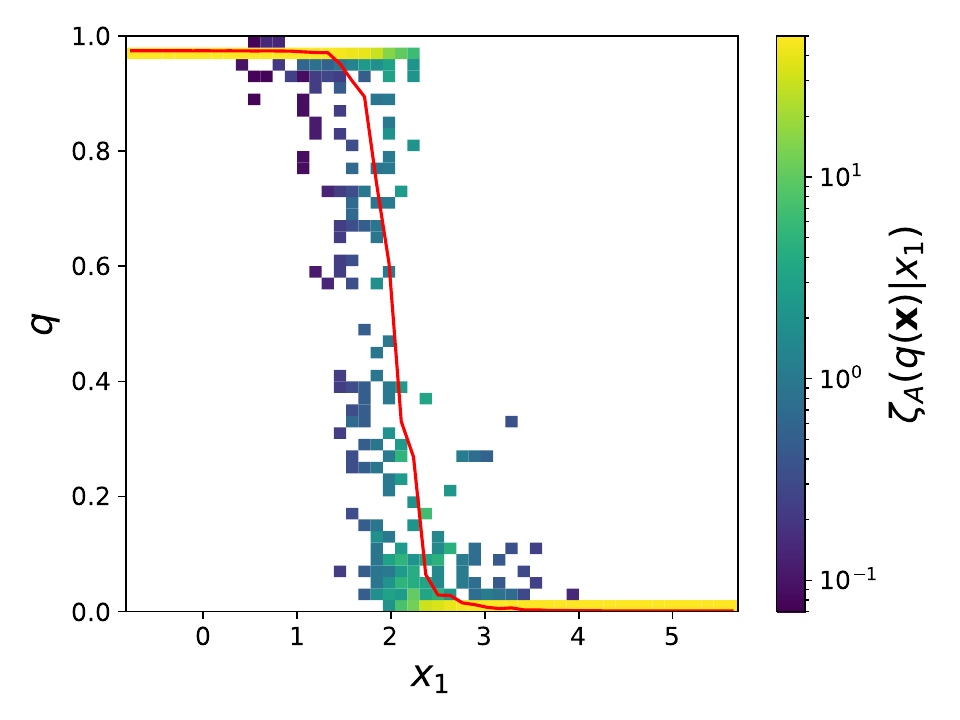} }\quad
\subfloat[][\emph{$15$ reactive trajectories.}\label{CommittorAnalogue_15}]
	{\includegraphics[width=.48\textwidth]{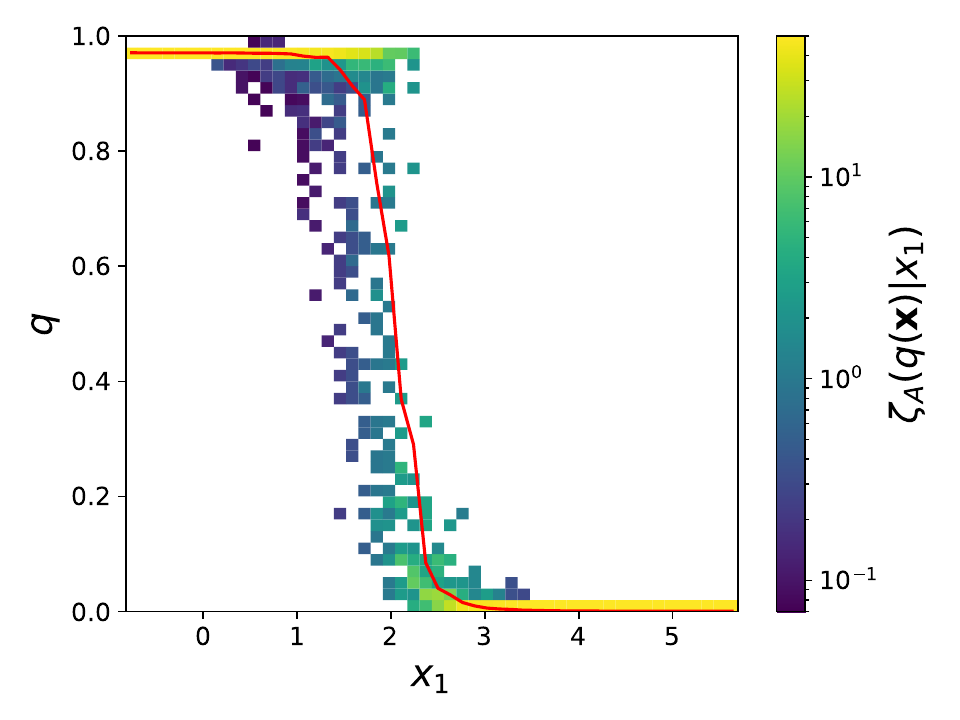}}\\
\caption{
Color maps of the conditional distribution $\zeta_A(q|x_1)$ for the committor function of the Charney--DeVore model estimated using the analogue method $\hat{q}^A(\mathbf{x})$ with learning dataset containing 2 (a) and 15 (b) reactive trajectories. The conditional average $\langle q_A\rangle_{x_1}$ as a function of $x_1$ is shown as a red line.
}
\label{fig:CommittorGeneratorRegressor_CDV}
\end{figure*}
We represent the conditional distributions $\zeta_A(q|x_1)$ in figure~\ref{fig:CommittorGeneratorRegressor_CDV}, along with the conditional average $\langle q_A\rangle_{x_1}=\int q\zeta_A(q|x_1)\,{\rm d}q$.
The subscript ${}_A$ indicates that they are conditional distributions for the estimated committors $\hat{q}_A$, using the analogue Markov chain. The results are shown for two learning experiments, using time series displaying respectively 2 and 15 reactive trajectories.
By comparing figure~\ref{fig:CommittorReferenceN100} and figure~\ref{fig:CommittorGeneratorRegressor_CDV}, it can be noted that the conditional distributions $\zeta_A(q|x_1)$ provided by the analogue method have the same qualitative structure as the conditional distribution of the reference committor, with values concentrated close to 0 and 1 in the vicinity of sets $\mathcal{A}$ and $\mathcal{B}$, and a sharp transition region in between.
However, the distributions are much more concentrated around the two set $\mathcal{A}$ and $\mathcal{B}$ than the reference one.
This is probably because the phase space has not been explored sufficiently and therefore the analogues of points lying outside the hyperballs defining the sets are instead inside $\mathcal{A}$ and $\mathcal{B}$.
Similarly, the transition region is narrower.
The estimates obtained with the two datasets of different lengths are qualitatively very similar (see Figs.~\ref{CommittorAnalogue_2} and~\ref{CommittorAnalogue_15}), even if the distribution using 15 reactive trajectories (figure~\ref{CommittorAnalogue_15}) exhibits slightly more spread close to attractor $\mathcal{B}$ and a seemingly broader transition region.

We now compare the performances of the two data-based methods (the analogue method and the direct estimator~\eref{Eq:CommittorTrajectory}) as the amount of data varies using the same procedure as in Sec.~\ref{sec:example2D}.
The error associated to an estimate of the committor is given by the non-constant term of the Brier score~\eref{Eq:TotalBrier}, i.e.
\begin{equation}
\| q-\hat{q}\|^2_{\rho_s}=\int{\rm d}\mathbf{x} {(q(\mathbf{x})-\hat{q}(\mathbf{x}))}^2\rho_s(\mathbf{x}) \approx \frac{1}{N_{p}}\sum_{i=1}^{N_{p}}{(q(\mathbf{x}_i)-\hat{q}(\mathbf{x}_i))}^2,
\label{eq:error_cdv}
\end{equation}
where $q$ is the reference committor, $\hat{q}$ its approximation and $\rho_s$ is the invariant measure.
Note that here, we are directly comparing the committor functions $q$ and not the distributions $\zeta(q|x_1)$.

For each dataset size, we repeat the computation 10 times using different realizations of the trajectory.
The error is computed as the empirical average over those realizations and the error bar corresponds to the standard deviation computed over the different experiments.
\begin{figure}[tp]
  \centering
  \includegraphics[width=.5\textwidth]{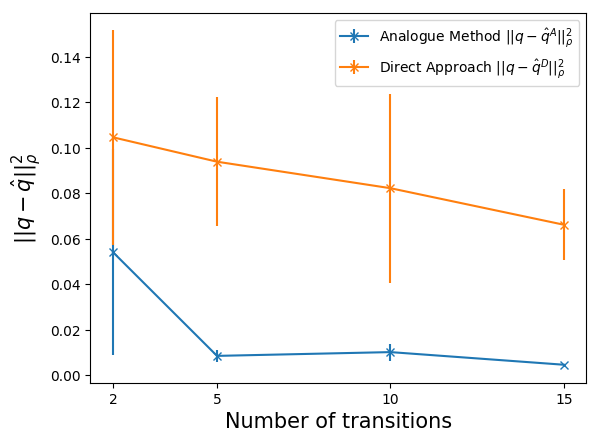}
  \caption{Error for the analogue and direct estimators of the committor function of the Charney--DeVore model as function of the number of reactive trajectories in the learning dataset.}
  \label{fig:Error_Committor_CDV}
\end{figure}
These results are shown, as a function of the size of the dataset upon which the analogue Markov chain is built, in figure~\ref{fig:Error_Committor_CDV}.
The estimates of the committor function provided by the analogue method are more precise than those obtained by the direct approach, regardless of the length of the dataset.
While the error associated to the direct approach decreases linearly with the number of transitions in the dataset, the error for the analogue method decreases faster for small datasets and reaches a plateau after about 5 reactive trajectories.
In this latter regime, the error with the analogue method is roughly an order of magnitude smaller than with the direct approach.
It is remarkable that although the system has a larger number of degrees of freedom than the 2D system studied in Sec.~\ref{sec:example2D}, the error associated to the analogue method is likewise small for datasets containing a relatively small number of trajectories (about 5).
In other words, this suggests that the analogue method does not require larger datasets for estimating committor functions when the dimension of space increases.
Here, we have only verified this statement for a moderate increase of the number of degrees of freedom, but we hope that it remains true in higher dimension, as the relevant data should remain close to the transition paths where most of the information carried by the committor function is contained.

\section{Using the learned committor function in Adaptive Multilevel Splitting}\label{sec:ams}

In Sec.~\ref{sec:analogue}, we estimated the committor function with the analogue method.
We will now illustrate how this approximated committor can be used in a rare event simulation, using the Adaptive Multilevel Splitting (AMS) algorithm.
This algorithm relies on a function used to select the trajectories leading to the rarest events, called the \emph{score function}.
The committor function is known to be the optimal score function, but it is generally not known exactly.
We will show that using the estimated committor as a score function has two advantages.
First, it provides a version of AMS where the user does not need to explicitly prescribe the score function.
This is very useful in practice when we have little knowledge of the dynamics beyond the presence of the two attractors $\mathcal{A}$ and $\mathcal{B}$.
In addition, it can improve the precision of the quantities computed with AMS, compared to user-defined score functions.
Indeed, it approximates the true committor, which leads to minimal errors on estimates.

\subsection{The Adaptive Multilevel Splitting algorithm and the quality of score functions}\label{sec:amsalgo}

Adaptive Multilevel Splitting is a \emph{splitting} method designed to estimate the probability of rare events, inspired by the pioneering works of Kahn and Harris~\cite{kahn1951estimation} and Rosenbluth and Rosenbluth~\cite{rosenbluth1955monte}.
It has been proposed by C\'erou \& Guyader~\cite{cerou2007adaptive}, as an improvement over Multilevel Splitting (see~\cite{glasserman1998large} for instance).
Many variants have been developed since, and the algorithm has been applied in a variety of contexts~\cite{rolland2018extremely, bouchet2019rare, lopes2019analysis, lestang2020numerical}.
The description of the algorithm given here follows the presentation of Lestang \emph{et al.}~\cite{lestang2018computing}.
See the review article by C\'erou, Guyader \& Rousset~\cite{cerou2019adaptive} for a recent overview of the method and its applications.

For definiteness, we consider a continuous time Markov process $X_t$ in the phase space $\mathcal{X}$.
Let us define two regions $\mathcal{A}$ and $\mathcal{B}$ in phase space.
We again seek to estimate the probability $\alpha = \mathbb{P}[T_\mathcal{B} < T_\mathcal{A}]$, where $T_\mathcal{D} = \inf \{ t>0, X_t \in \mathcal{D} \; {\rm with} \; X_0 \in \mathcal{C} \}$ is the first hitting time of the set $\mathcal{D}$, starting from a set $\mathcal{C}$.
The set $\mathcal{C}$ encloses the set $\mathcal{A}$.
We also wish to compute the corresponding realizations of the dynamics.

The AMS algorithm computes these quantities iteratively.
For this matter, the algorithm uses a \emph{score function} $\phi$, (sometimes termed \emph{reaction coordinate}) a map from the phase space $\mathcal{X}$ to $\mathbb{R}$.
Ideally, the score function is bounded from below by $0$ and from above by $1$, vanishes identically on $\mathcal{A}$ and is identically equal to 1 on $\mathcal{B}$.
Our aim is to compare the efficiency of different score functions.

In order to run the algorithm, we first need to sample initial conditions according to the invariant measure restricted to the set $\mathcal{C}$.
In practice, we sample these initial conditions on $\mathcal{C}$ by sampling long trajectories in the basin of attraction of $\mathcal{A}$.
Then the algorithm is initialized by sampling $N$ independent trajectories, with initial conditions on the set $\mathcal{C}$ and run until they reach either the set $\mathcal{A}$ or the set $\mathcal{B}$.
Let us denote by ${\{\mathbf{x}_n^{(0)}(t)\}}_{1 \leq n \leq N}$ the initial ensemble of trajectories, where the subscript denotes the index of the trajectory in the ensemble and the superscript denotes the iteration of the algorithm.
We associate a weight $w_0=1$ to those trajectories.

\begin{figure*}[tp]
  \centering
  \includegraphics[scale=.5]{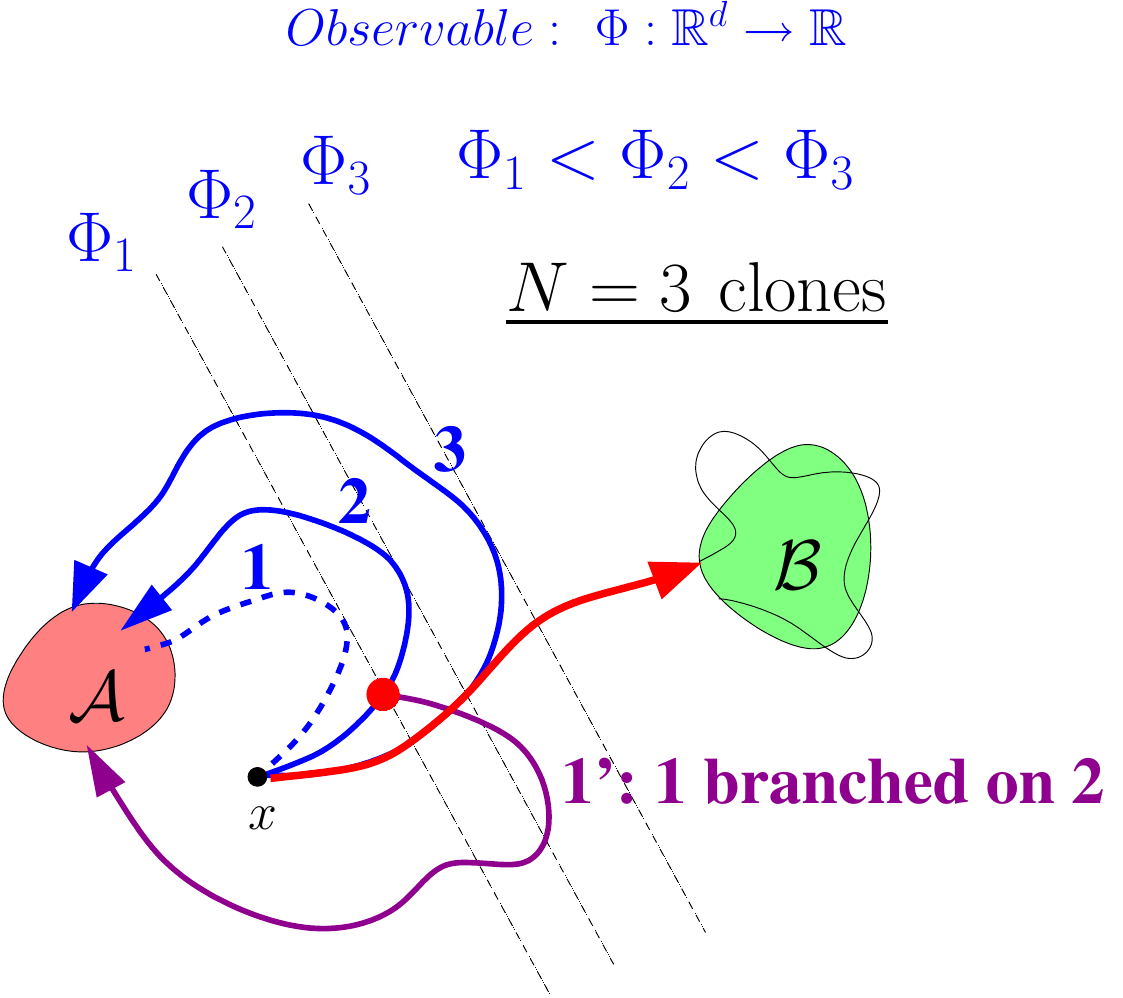}
  \caption{
    Sketch illustrating two iterations of AMS in a simplified example with $3$ clones (Figure originally made for~\cite{simonnet2016combinatorial}), in order to compute trajectories going from set $\mathcal{A}$ to set $\mathcal{B}$.
    Trajectory $1$ (dashed line) has the smallest excursion out of $\mathcal{A}$ as measured by the score function $\Phi$.
    It is removed and branched on another trajectory (in that case trajectory 2, leading to the purple line).
    In the successive iteration, trajectory 2 has the smallest score function and is branched on trajectory 3 (leading to the red line).
  }
  \label{fig:Illustration_AMS}
\end{figure*}
At each iteration $j \geq 1$, we apply the following \emph{selection and mutation} steps, which are schematically illustrated in figure~\ref{fig:Illustration_AMS}:
\begin{itemize}
\item We compute the score of each trajectory in the ensemble at iteration $j-1$: $\Phi_n^{(j)}= \sup_t \phi(t, \mathbf{x}_n^{(j-1)}(t))$.
\item We determine the trajectories which have the lowest score: $\Phi_j^\star= \min_{1\leq n \leq N} \Phi_n^{(j)}$ and we set $n_{j,1}^\star,\ldots,n_{j,\ell_j}^\star$ the indices such that $\Phi_{n_{j,1}^\star}^{(j)} = \cdots = \Phi_{n_{j,\ell_j}^\star}^{(j)} = \Phi_j^\star$. One can have $\ell_j>1$ in some iterations. If $\ell_j=N$ and not all the trajectories have reached $\mathcal{B}$, the algorithm stops: it leads to an \emph{extinction}.
\item We mutate each trajectory $\mathbf{x}_{n_{j,\ell}^\star}^{(j-1)}$ ($1 \leq \ell \leq \ell_j$): for each of them, we choose a trajectory $\mathbf{x}_{n_\ell}^{(j-1)}$ ($n_{\ell} \neq n_{j,1},\ldots n_{j,\ell_j}$) drawn randomly among the $N-\ell_j$ remaining trajectories. We determine the smallest time $t$ such that $\phi(t,\mathbf{x}_{n_\ell}^{(j-1)}(t))>\Phi_j^\star$, denoted by $t_{j,\ell}$. The new trajectory $\mathbf{x}_{n_{j,\ell}^\star}^{(j)}$ is set by copying the trajectory $\mathbf{x}_{n_\ell}^{(j-1)}$ from $t_0$ to $t_{j,\ell}$, and simulating the trajectory with a new independent realisation of the noise, starting from time $t_{j,\ell}$, until it hits either the set $\mathcal{A}$ or the set $\mathcal{B}$.
\item Trajectories with higher scores are not modified at this step: $\mathbf{x}_n^{(j)}=\mathbf{x}_n^{(j-1)}$ for $n \neq n_{j,1}^\star,\ldots,n_{j,\ell}^\star$.
\item We compute the weight of iteration $j$: $w_j=\left( 1 - \frac{\ell_j}{N}\right)w_{j-1}$.
\end{itemize}
The algorithm is iterated until all the trajectories reach the set $\mathcal{B}$.
The number of iterations $J$ is a random number.
This leads to an estimator $\hat{\alpha}$ for the transition probability $\alpha$:
\begin{equation}
  \hat{\alpha} = w_J= \prod_{j=0}^J\left(1-\frac{\ell_j}{N}\right)\,.\label{esta}
\end{equation}
This estimator is a random variable, with one value obtained for each realization of the algorithm.
We perform $M$ independent realizations of the algorithm and compute the statistics of $\hat{\alpha}$: the empirical average and variance of $\hat{\alpha}$.

The mathematical properties of this estimator have been extensively studied~\cite{cerou2007adaptive, guyader2011simulation,rolland2015statistical, brehier2015large, brehier2015analysis, brehier2016unbiasedness, brehier2016central, simonnet2016combinatorial}.
The key property is that, for any $N$ and score function $\phi$, it is an unbiased estimator~\cite{brehier2016unbiasedness, brehier2016central} with a finite variance.
The standard deviation $\sigma_\alpha(N)$, the square root of the variance, depends on $N$ and on the score function.
The optimal score function, with the lowest variance, is the committor function.

More precise results exist asymptotically for large $N$.
It is then proven~\cite{cerou2019asymptotic} that the standard deviation scales like $\frac{1}{\sqrt{N}}$ asymptotically $\sigma_\alpha(N) \sim\frac{G(\phi)}{\sqrt{N}}$.
Moreover, when the score function is the committor function, $G$ is minimal, and the standard deviation scales like the ideal standard deviation
\begin{equation}
\sigma_{\rm id}=\frac{\alpha \sqrt{|\log(\alpha)|}}{\sqrt{N}}\,.\label{esid}
\end{equation}
In many cases, an asymptotic scaling is observed in practice when the number of clones is larger than $100$  (see for instance~\cite{rolland2018extremely}, figure~14 (c)).
The computation of the empirical standard deviation of $\alpha$, given by
\begin{equation}
\sigma_\alpha(N,M)=\sqrt{\frac{1}{M}\sum_{m=1}^M\left( \hat{\alpha}_m^2\right)-{\left(\frac{1}{M}\sum_{m=1}^M \hat{\alpha}_m \right)}^2}\,,
\end{equation}
and its comparison to the ideal standard deviation $\sigma_{\rm id}$ has often been used as an \emph{a posteriori} test of the quality of the score function and how close it is to the committor~\cite{rolland2015statistical,brehier2019new,rolland2016computing}.

Although the estimator is actually unbiased ($\mathbb{E}[\hat{\alpha}]=\alpha$), in numerical uses of AMS, it is often observed that $\hat{\alpha}$ underestimates $\alpha$ in the large majority of the $M$ realizations of the algorithm.
These underestimates are such that the average $\langle\hat{\alpha}\rangle_M=\frac{1}{M}\sum_{m=1}^M\hat{\alpha}_m$ over $M$ realizations is most of the time strictly smaller than $\alpha$ although the average is $\alpha$ ($\mathbb{E}[\langle\hat{\alpha}\rangle_M] = \alpha$).
This phenomenon is called an \emph{apparent bias}.
We note that a similar observation is made in the context of fixed Multilevel Splitting~\cite{glasserman1998large} and Importance Sampling~\cite{devetsikiotis1993statistical}.
In these contexts, it can be demonstrated that $\frac{1}{M}\sum_{m=1}^M\hat{\alpha}_m$ underestimates $\alpha$ with a probability that goes to $1$ as parameters like $\epsilon$, which control the rareness of the event, go to zero~\cite{glasserman1998large}.
This happens if the score function yielding the levels of Multilevel Splitting is not adapted.
As a consequence, the observed sample mean of $\hat{\alpha}$ will be strictly smaller than $\alpha$ unless an out of reach number of realizations of AMS is performed.
It has been conjectured~\cite{brehier2016unbiasedness} that the observed apparent bias phenomenon could be explained for the AMS by analogy with the studies for fixed multilevel splitting.

The apparent bias, measured through the difference $\alpha-\langle \hat{\alpha}\rangle_M$, decreases like $\frac{1}{N}$ as the number of clones $N$ is increased.
However, it has been observed that for some cases, the apparent bias seems to reach a plateau for extremely large values of $N$~\cite{rolland2015statistical}.
We will see similar behavior in the following.
In these situations, it is observed that this apparent bias is minimal when the score function is the committor function~\cite{rolland2015statistical}.
As this apparent bias is a very important practical problem, we will use the magnitude of this apparent bias as a measure of the quality of the score function.

We have seen that the committor function is the best score function for the AMS algorithm, and explained that the computation of the empirical standard deviation and of the apparent bias are two ways to quantify the quality of a score function.
We can also test the AMS computations by comparing the computation of other observables.
For instance, we will compute the transition path duration, denoted $\tau$.
This physical quantity has proven to be a good indicator of whether AMS was correctly sampling transition paths~\cite{rolland2015statistical}.

To have an unbiased estimate of $\alpha$ and $\tau$ and validate the output of AMS computations, we perform a large number of Direct Numerical Simulations (DNS) of reactive trajectories.
These DNS start like AMS computations with initial conditions on $\mathcal{C}$, we let them evolve until they reach either $\mathcal{A}$ or $\mathcal{B}$.
The proportion of DNS that reach $\mathcal{B}$ before $\mathcal{A}$ yields a direct estimate of $\alpha$.
We also perform an estimate of $\tau$ by averaging the duration of trajectories that reach $\mathcal{B}$ before $\mathcal{A}$.

The estimate of a quantity by AMS is deemed to be precise enough when the 95\% confidence intervals of this estimate performed by AMS and by DNS overlap~\cite{brehier2016unbiasedness}.
These confidence intervals are constructed by noting that we look at the sum of independent random variables of finite variance.
They therefore follow a central limit theorem and the sample mean of $\hat{\alpha}$ has a gaussian distribution.
The confidence interval is then given by $\langle \alpha\rangle_M\pm 1.96\sigma_\alpha(N,M)$, with the empirical standard deviation $\sigma_\alpha(N,M)$.
Similar confidence intervals are constructed for $\alpha$ and $\tau$ for both AMS and DNS results.

\subsection{The learned committor function}\label{sec:datascorefunction}

Our goal is to investigate the performance of a score function relying on a data-based estimate of the committor function, using the analogue method presented in Sec.~\ref{sec:analogue}.
As mentioned above, this method only provides an estimate on the points initially present in the dataset.
To extend the score function to the whole phase space, we proceed as explained in Sec.~\ref{sec:extendingcommittor}, with a nearest-neighbor method using an exponential kernel with width $\omega=0.1$.
Here, a small number of neighbors $K=10$ is used for efficient computations of the score function.
Indeed, for each computation, a search through neighbors must be performed.
For a given training dataset, this method defines a score function, which we denote $\phi_{\rm dat}$.

The use of a kernel is justified by the need to reduce regions of constant score function.
Indeed $\mathbb{R}^D$ is divided in finite subvolumes where any point $\mathbf{y}$ has the same neighbours $\{\mathbf{x}_j\}_{1\le j\le K}$ and thus would have the same score function without a kernel.
The use of the kernel ensures that sets of constant $\phi_{\rm dat}$ are hypersurfaces and not hypervolumes and that much fewer $\mathbf{y}$ have the same values of $\phi$ at each stage of the algorithm.
Practice shows that this leads to more efficient branching by limiting the number of clones suppressed at each stage of the algorithm and the risk of extinction.

We will test the use of the analogue based estimate of the committor as a score function for the AMS computations for the two systems presented in Sec.~\ref{sec:committorapplications}: the 2D three-well system (Sec.~\ref{sec:amsthreewells}) and the Charney--DeVore model (Sec.~\ref{sec:amscdv}).
The test of the learned committor function will be twofold.
First, we will consider a score function learned on a dataset displaying a large number of transitions and study the quality of the result as a function of the number of clones $N$.
This will allow us to discuss the phenomenon of apparent bias and how the learned committor function deals with it.
The second aim will be to study the required size of the dataset to have good results with the AMS algorithm.
This question is critical for complex systems for which data will be scarce because of computation costs.
To address this question, we will then perform AMS computations with a fixed large number of clones $N=1000$ and datasets of increasing size (measured in number of recorded transitions).

\subsection{AMS study for the two dimensional three well model}\label{sec:amsthreewells}

In this subsection, we work on the dynamics of the two-dimensional three-well model presented in Sec.~\ref{sec:example2D}~\eref{eq:Dynamics}.
The sets $\mathcal{A}$ and $\mathcal{B}$ as well as the noise variance $\epsilon$ are defined as in Sec.~\ref{sec:example2D}.

\subsubsection{Efficiency of the AMS algorithm with the learned committor function for large $N$ for the three-well model}\label{sec:apparentbias3wells}

We first study the efficiency of the AMS algorithm when using the learned committor function $\phi_{\rm dat}$ as a score function, when the number of clones is increased with a fixed data set length.

The time series which is used to compute this score function has $N_t=1.4\cdot 10^5$ datapoints (effectively 1400 time units long) and displays 21 transitions.
The results for the AMS algorithm with this score function will be compared either to DNS computations, or to AMS computations with two explicit user-defined score functions: $\phi_{\rm lin}(\mathbf{x})=\frac{x+1}{2}$ and $\phi_{\rm norm}(\mathbf{x})=\frac{\sqrt{{(x+1)}^2+\frac12 y^2}}{2}$.
The performances of these score functions have been studied in detail in the literature~\cite{rolland2015statistical}.
The sample means are performed using at least $M=6000$ realisations of AMS (with $\phi_{\rm dat}$ and $N=1000$) and up to $M=250000$ (with $\phi_{\rm lin}$ and $N=10$).

It has been observed in the literature that if one use the first score function $\phi_{\rm lin}$, the probability $\alpha$ can be gravely underestimated if one uses a number of realisations $M$ too small for a given number of clones. In practice, at fixed $M$ and $N$ the sample mean $\frac{1}{M}\sum_{m}\hat{\alpha}_m$ goes to zero with the noise amplitude $\epsilon$ with a larger rate than the one for $\alpha$.
 The sample mean of durations of reactive trajectories would also be strongly underestimated.
 With that score function, AMS computations wrongly overselect trajectories going through the bottom channel (where paths remain around $y\simeq 0$ and where they cross the highest potential difference, see figure~\ref{fig:Potential} and~\cite{rolland2015statistical}, Fig. 7)a)).
As a consequence the most probable value of $\hat{\alpha}$ is much smaller than $\alpha$.
 It is conjectured that these severe underestimations are a consequence of the apparent bias phenomenon for AMS~\cite{brehier2016unbiasedness}.
The correct estimate is recovered asymptotically in the limit $M\rightarrow \infty$.
For this to happen, the histograms of $\hat{\alpha}$ develop large power law tails toward very large values of $\hat{\alpha}$.
This phenomenon already starts to be at play for our case of $\epsilon=0.5$.
What is at stake here is to propose a score function that ensures the most precise estimate while requiring as few clones and realisations as possible.
We will show that our learned committor score function does lead to this improvement.

\begin{figure*}
\centerline{\includegraphics[width=0.33\linewidth]{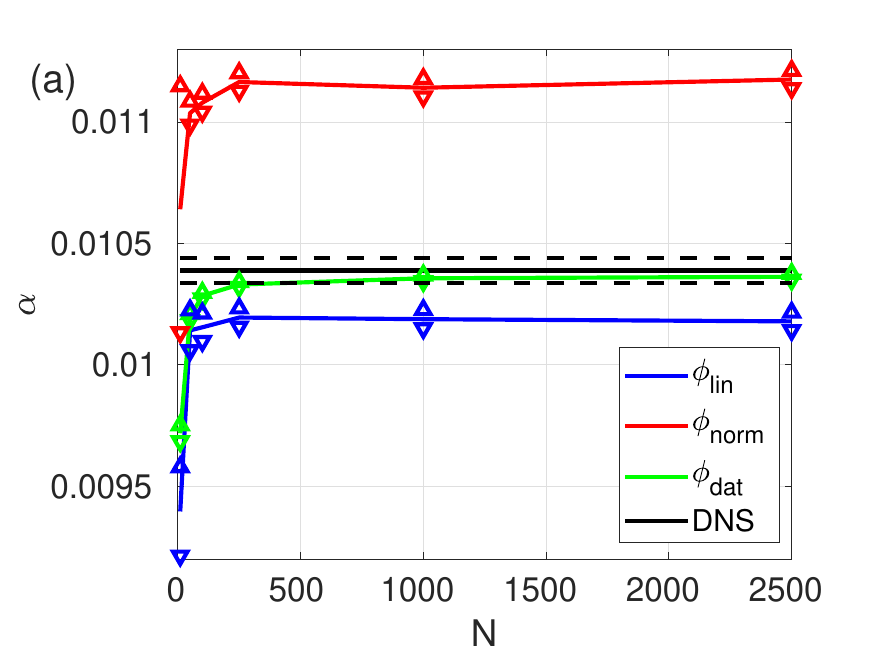}\includegraphics[width=0.33\linewidth]{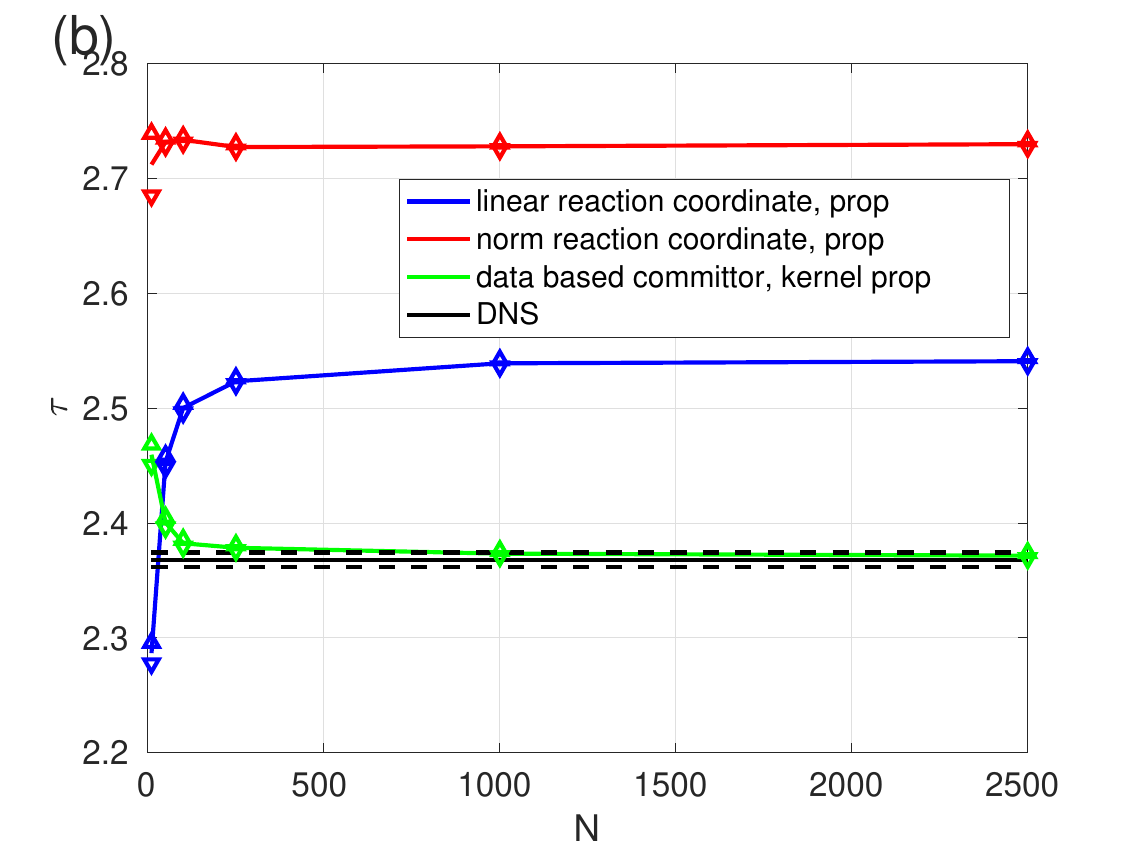}\includegraphics[width=0.33\linewidth]{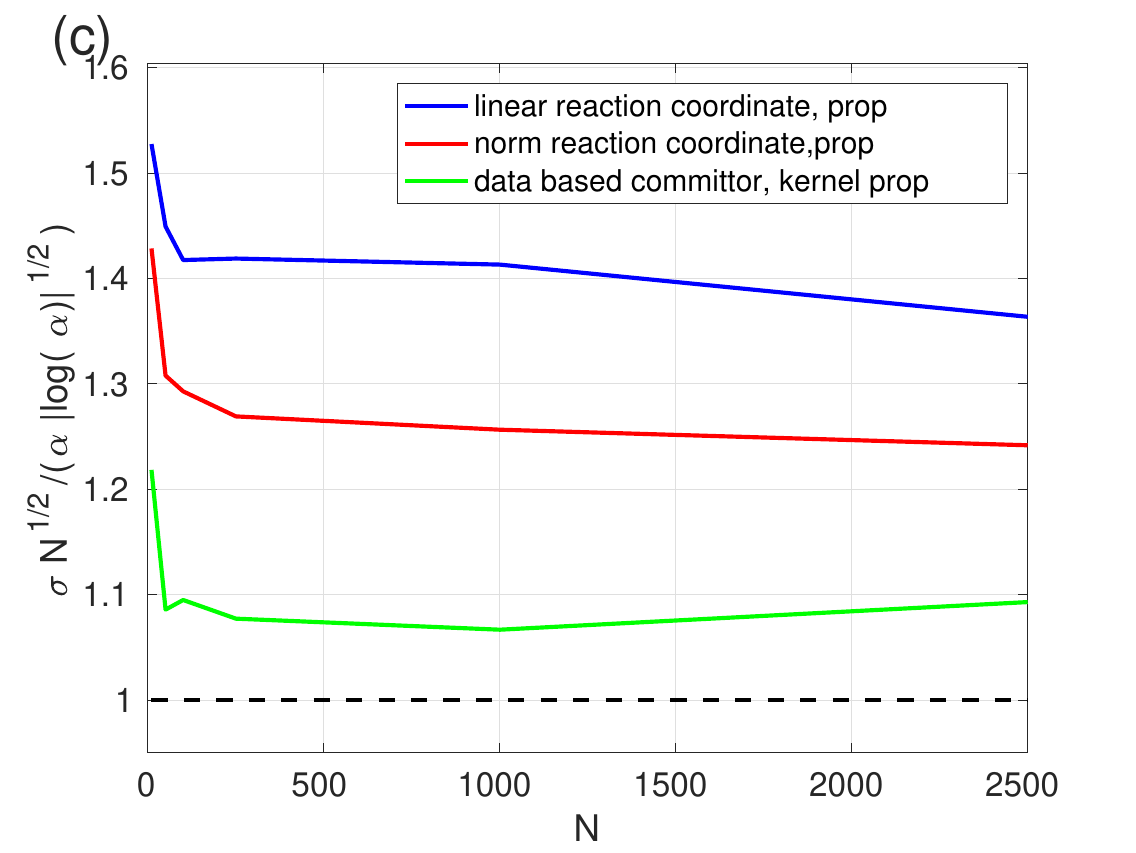}}
\caption{Efficiency of the AMS algorithm with the learned committor function for fixed and large dataset, for the three-well problem. Comparison of the estimated (a) Transition probability $\langle \alpha\rangle$, (b) Duration of reactive trajectories $\langle \tau\rangle$, and (c) Rescaled standard deviation $\sigma$, as a function of the number of clones $N$. For each plot the black curve is the reference: either the DNS (a) and (b), or the optimal value 1 (c). The dashed black lines are the 95\% confidence interval for the DNS. The color curves have been computed using the AMS, with respectively the learned committor function (green), the linear score function (blue) and the quadratic score function (red). The red and blue curves clearly illustrate the apparent bias phenomenon. The learned committor function gives excellent results, suppressing the apparent bias and giving smaller, close to optimal, empirical standard deviation.}
\label{eval_temp}
\end{figure*}
In figure~\ref{eval_temp} (a), we first show the transition probability $\langle \alpha\rangle_M$ as a function of the number of clones $N$ used in AMS computations, using the three score functions $\phi_{\rm dat}$, $\phi_{\rm lin}$ and $\phi_{\rm norm}$, and computed by means of DNS.
Error bars show the 95\% interval of confidence.
One can first note that for the AMS computations, $\langle\alpha\rangle_M(N)$ is within $1\%$ of its asymptotic value if the number of clones used is larger than $N=100$.
As noted in Sec~\ref{sec:amsalgo} $\langle\alpha\rangle_M$ grows with $N$ toward this asymptotic value.
The confidence intervals of the probability for AMS and DNS computations do not overlap when we use the norm score function $\phi_{\rm norm}$: the asymptotic value of $\langle\alpha\rangle_M$ overestimates $\alpha$.
With the linear score function $\phi_{\rm lin}$, the asymptotic value of $\langle \alpha\rangle_M$ in turn underestimates $\alpha$: this is a possible consequence of the apparent bias phenomenon~\cite{brehier2016unbiasedness}.
These results are in agreement with previous studies~\cite{rolland2015statistical}, which have related these biases to errors in the relative sampling of transition paths.
For instance, the linear score function selects preferentially trajectories going through the bottom channel (where paths remain around $y\simeq 0$ and where they cross the highest potential difference, see figure~\ref{fig:Potential} and~\cite{rolland2015statistical}, Fig. 7)a)), leading to the bias.
By contrast, if the learned committor function $\phi_{\rm dat}$ is used, the confidence intervals overlap as soon as $N\ge 250$, thus indicating that no bias can be detected in this estimate of $\alpha$.

The results for $\tau$, the average duration of reactive trajectories, are qualitatively similar: figure~\ref{eval_temp} (b) shows $\langle \tau\rangle_M$ as a function of the number of clones used in AMS for the three score functions, compared to a reference DNS calculation.
Error bars are again given by the 95\% confidence interval.
We first note that for the data-based and linear score functions $\langle\tau\rangle_M$ converges toward its asymptotic value to within 1\% for $N \ge 250$; for the norm score function it is within that interval for all values of $N$.
The 95\% confidence interval of the AMS estimate with learned committor function and of DNS overlap if more than $1000$ clones are used (a larger number than for the transition probability $\alpha$; note however that the confidence intervals are narrower for $\tau$ than for $\alpha$).
Conversely, both the linear and the norm score functions lead to overestimates of $\tau$.

Finally, in figure~\ref{eval_temp} (c), we consider the rescaled standard deviation $\sigma_{\alpha}(N,M)/\sigma_{\rm id}(N)$ of the estimator of $\alpha$ as a function of the number of clones  $N$.
Here the absolute reference is the unit value, obtained for the optimal score function, the exact committor.
We first note that for all score functions, the rescaled standard deviation reaches a plateau if the number of clones is larger than $N=100$.
The value of this plateau is largest when we use the linear score function, with $\sigma_{\alpha}(N,M)=1.4\pm 0.02$.
It is somewhat smaller for the norm score function, with $\sigma_{\alpha}(N,M)=1.25\pm 0.02$. The best results are obtained for the learned committor function, with $\sigma_{\alpha}(N,M)= 1.12\pm 0.02$.
This again indicates that the computations performed using the learned committor function are the most precise, in that they come with the smallest statistical error, which is 10\% larger than the smallest error possible.

All things considered, we conclude that if we use a large dataset to learn an estimate committor function with the analogue Markov chain, and use it as a score function for AMS, the estimates of transition properties show no apparent biases and converge to their true value when the number of clones is increased.
The better precision of the results with the learned score function is also clearly visible for the lower statistical error measured by the empirical standard deviation of the estimate of $\alpha$. This level of precision is ensured if at least $N=1000$ clones are used: then the 95\% intervals of confidence of estimates by mean of DNS and by mean of AMS using the learned score function overlap.

\subsubsection{Efficiency of the AMS algorithm with the learned committor function as a function of the dataset length for the three-well model}\label{sec:datasize3wells}

In Sec.~\ref{sec:apparentbias3wells}, we used a large dataset with 21 transitions to accurately estimate the committor, before using it as a score function for the AMS.
Compared to analytically defined score functions, this suppressed the apparent bias phenomenon and reduced the statistical error.

However, for many complex systems with very costly computations, it might not always be affordable to use a long dataset to learn the committor function.
Moreover, in the initial stage of the study, one need to work with short datasets.
Hence, we now study how the results of AMS computations using the learned committor function depend on the size of the learning dataset. We also study how much the results change from one realization of the dataset to another, at fixed number of transitions. This study is first made for the 2D three-well model.

For this matter, we sample trajectories of increasing length that contain an increasing number of transitions, from 1 to 21.
For a number of transitions going from $2$ to $21$, we sample seven independent trajectories.
For each of these datasets, we estimate the committor with the analogue method (Sec.~\ref{sec:analogue}) and use it as a score function in AMS computations with $N=1000$ clones.
For all AMS computations using the same reaction coordinate ($\phi_{\rm lin}$, $\phi_{\rm norm}$ and $\phi_{\rm dat}$ with each dataset), outputs were averaged over $M=20000$ realisations.

\begin{figure*}
\centerline{\includegraphics[width=0.33\linewidth]{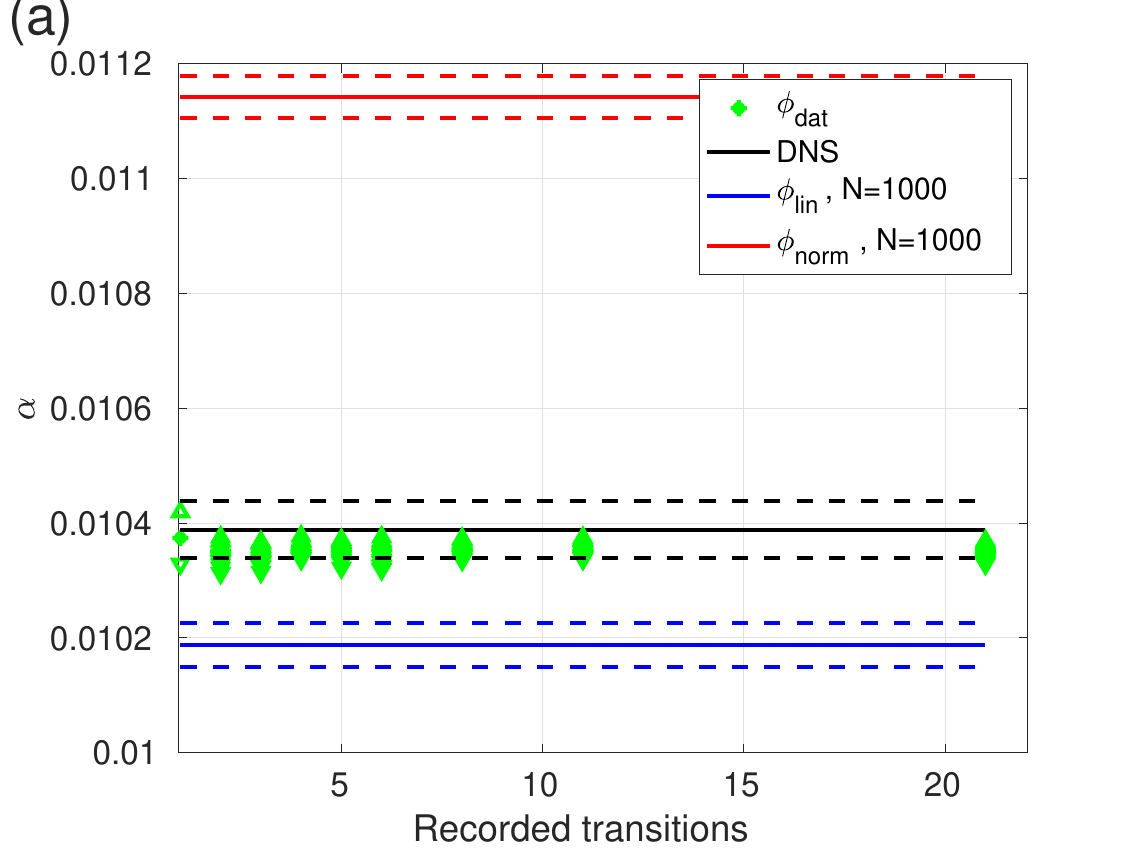}\includegraphics[width=0.33\linewidth]{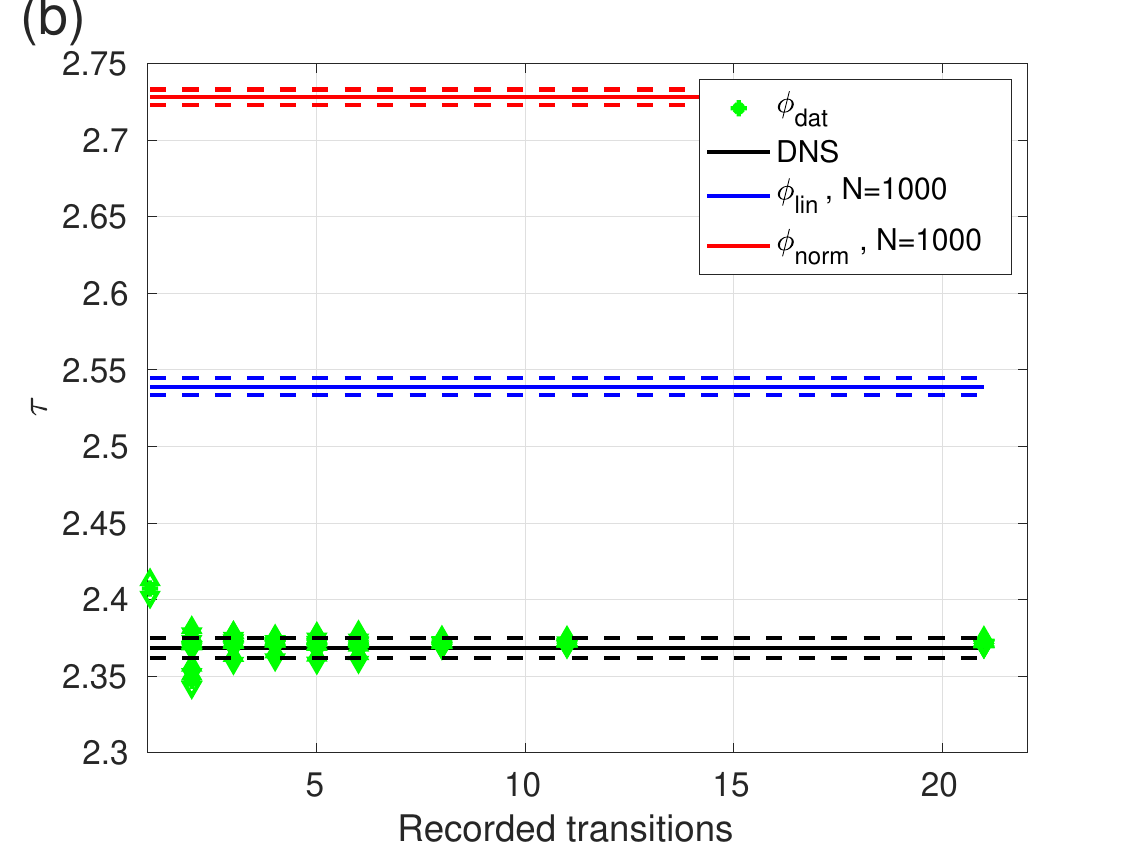}\includegraphics[width=0.33\linewidth]{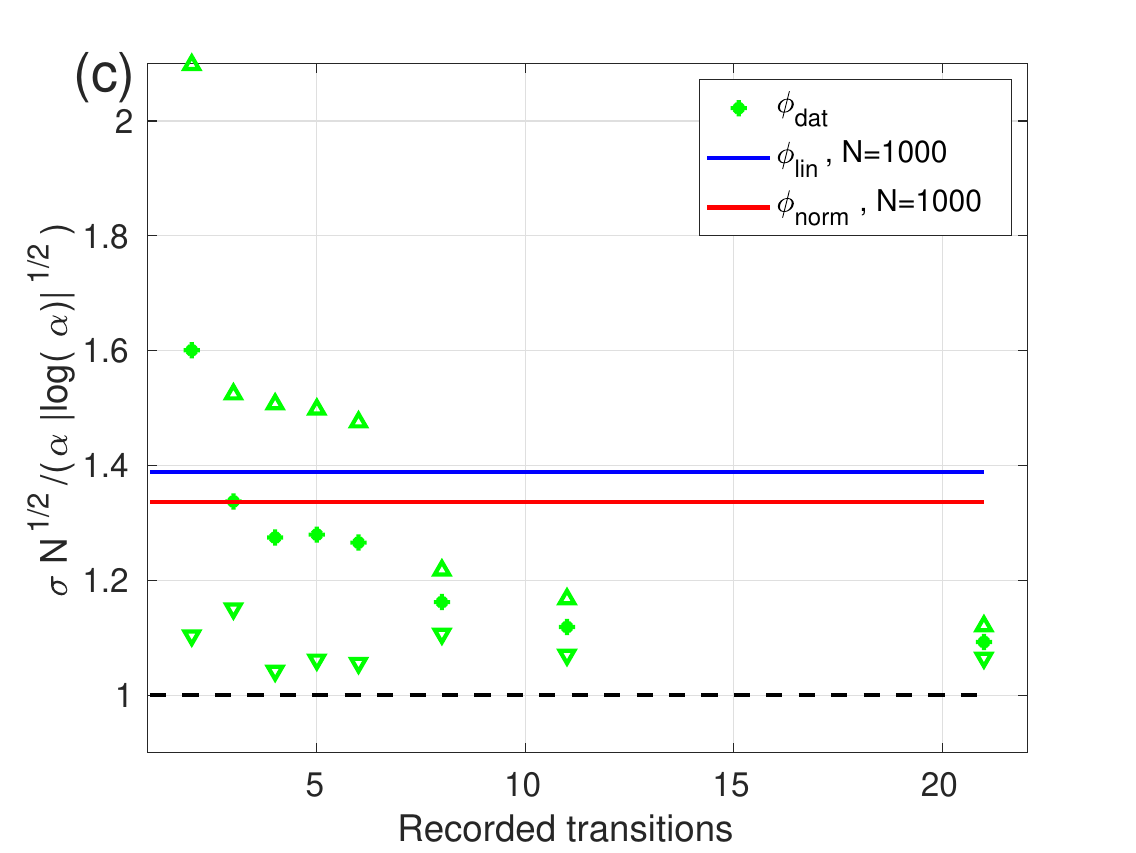}}
\caption{Efficiency of the AMS algorithm with the learned committor function, as a function of the dataset length, for the three-well problem. Comparison of the estimated (a) Transition probability $\langle \alpha\rangle$, (b) Duration of reactive trajectories $\langle \tau\rangle$, and (c) Rescaled standard deviation $\sigma$, for each case averaged over independent realizations of the score function. For each plot the black curve is the reference one, either the DNS (a and b) or the optimal value 1 (c). The dashed black lines are the 95\% confidence interval for the DNS. The color curves have been computed using the AMS, with respectively the learned committor function (green), the linear score function (blue) and the quadratic score function (red). The red and blue curves are constant values (they do not depend on the data set length) for comparison. The learned committor function gives much better results than the user-defined score functions, even for very small datasets. With datasets containing only a few transitions, two to five, the results are already excellent. However, for such small datasets, the quality of the score function varies much from one realization to another.}
\label{comm_datlen}
\end{figure*}
Figure~\ref{comm_datlen} (a,b) show the transition probability $\alpha$ and the average duration of reactive trajectories $\tau$ as a function of the number of sampled transitions in the dataset.
We place a point with the 95\% confidence interval as error bars of each AMS estimate using each distinct learned committor.
However, all the points are essentially superimposed: all the realizations of the score functions lead to the same results.
We note that even if we use a short dataset to learn the committor function, the estimates are very precise:
the intervals of confidence of the AMS algorithm and the DNS estimates overlap for all our datasets lengths, except for the shortest dataset (only one transition) for $\tau$.
Meanwhile the confidence interval of the estimate of $\tau$ by AMS using $\phi_{\rm lin}$ and $\phi_{\rm norm}$ do not overlap with that of DNS, which indicates the beginning of incorrect selection of trajectories.

This is confirmed in figure~\ref{comm_datlen} (c) by considering $\sigma$, the rescaled standard deviation of the estimate of $\alpha$.
In this plot, for each dataset length, we have computed the empirical average and variance of the rescaled standard deviation estimated with the different realizations of the score function.
The empirical average gives the points and the variance is used to construct the error bars.
This first shows that the rescaled standard deviation decreases as the number of transitions contained in the dataset increases, from 1.6 when the dataset contains only two transitions to almost 1.1 when the dataset contains 8 transitions or more.
This indicates that the statistical error obtained using the learned score function is systematically smaller than that obtained using the user-defined score function as soon as the datasets contain at least 8 transitions.

We also note that the fluctuations of the standard deviation between different dataset realizations decreases as the number of transitions contained in the dataset increases. If the dataset is short, no more than 6 transitions, one can obtain a score function that leads to better or worse results than analytically defined score functions with comparable probability. With a dataset with 3 transitions or more, the statistical error is most of time reduced when the score function is learned from datasets, compared to the case with user-defined score functions.

Finally, we stress that for very short datasets, with only a few transitions, even if the standard deviation on the estimate of $\alpha$ is of the same order for both user-defined and learned score function, the systematic apparent bias is much smaller with the learned committor function.

\subsection{Application to the Charney--DeVore model}\label{sec:amscdv}

We now perform the same tests for AMS computations using the learned committor function in the Charney--DeVore model~\eref{eq:CharneyDeVoreModel}.
We will compute transitions from zonal to blocked flows and use $r_Z=r_B=0.8$ to define the sets $\mathcal{A}$ and $\mathcal{B}$ as in Sec~\ref{sec:charneydevore}.

\subsubsection{Efficiency of the AMS algorithm with the learned committor function for large $N$ for the Charney--DeVore model}\label{sec:apparentbiasCDV}

We proceed as in Sec.~\ref{sec:apparentbias3wells}: we first learn the committor function $\phi_{\rm dat}$ (Sec.~\ref{sec:datascorefunction}) from a long trajectory, containing $N_t=3.4 \cdot 10^4$ data points and displaying 38 transitions.
The estimates will be compared to DNS results and to a simple linear score function $\phi_{x_1}=\frac{x_1-x_{1,Z}}{x_{1,B}-x_{1,Z}}$ with $x_{1,Z}=4.308$ and $x_{1,B}=0.709$ (see Figs.~\ref{ConvergenceToZonal} and~\ref{ConvergenceToBlocked}).
We will perform averages over at least $M=450$ realisations of AMS computations ($N=2500$ clones with the learned committor score function) and over up to $M=20000$ realisations ($N=10$ clones with the linear score function).

\begin{figure*}
\centerline{\includegraphics[width=0.33\linewidth]{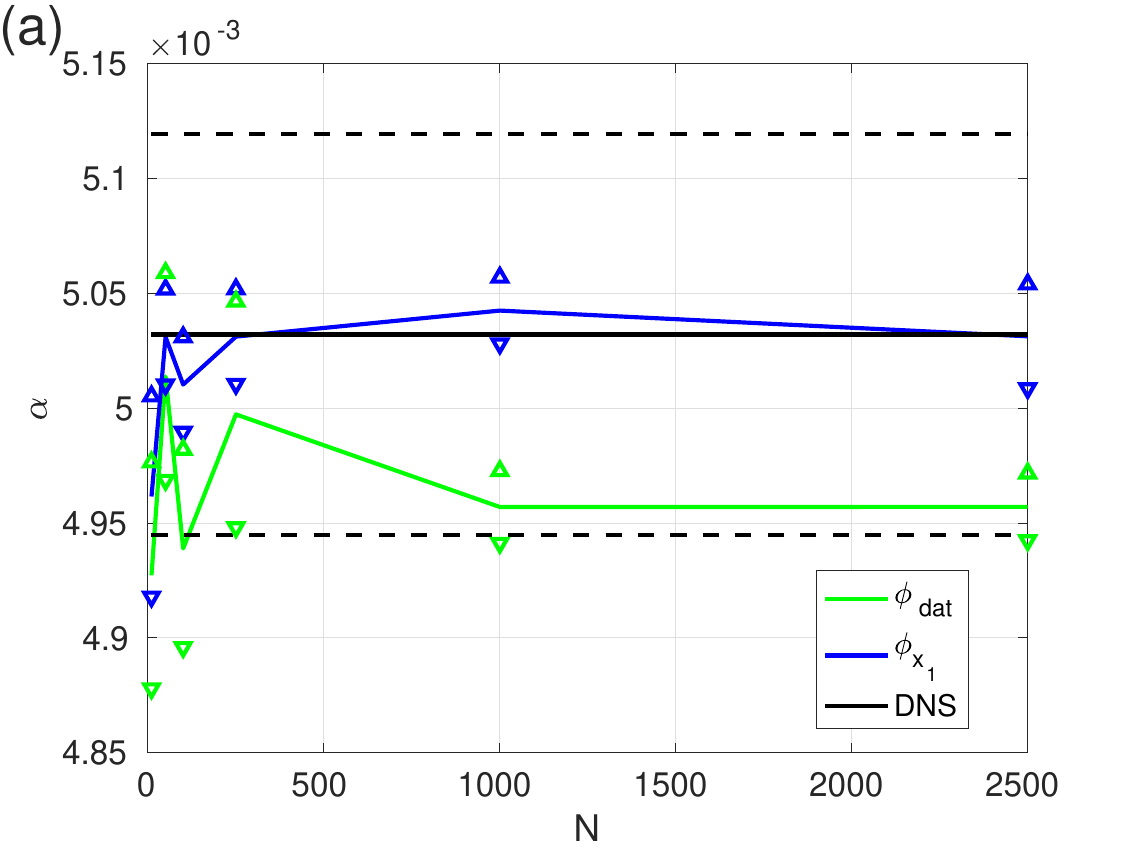}\includegraphics[width=0.33\linewidth]{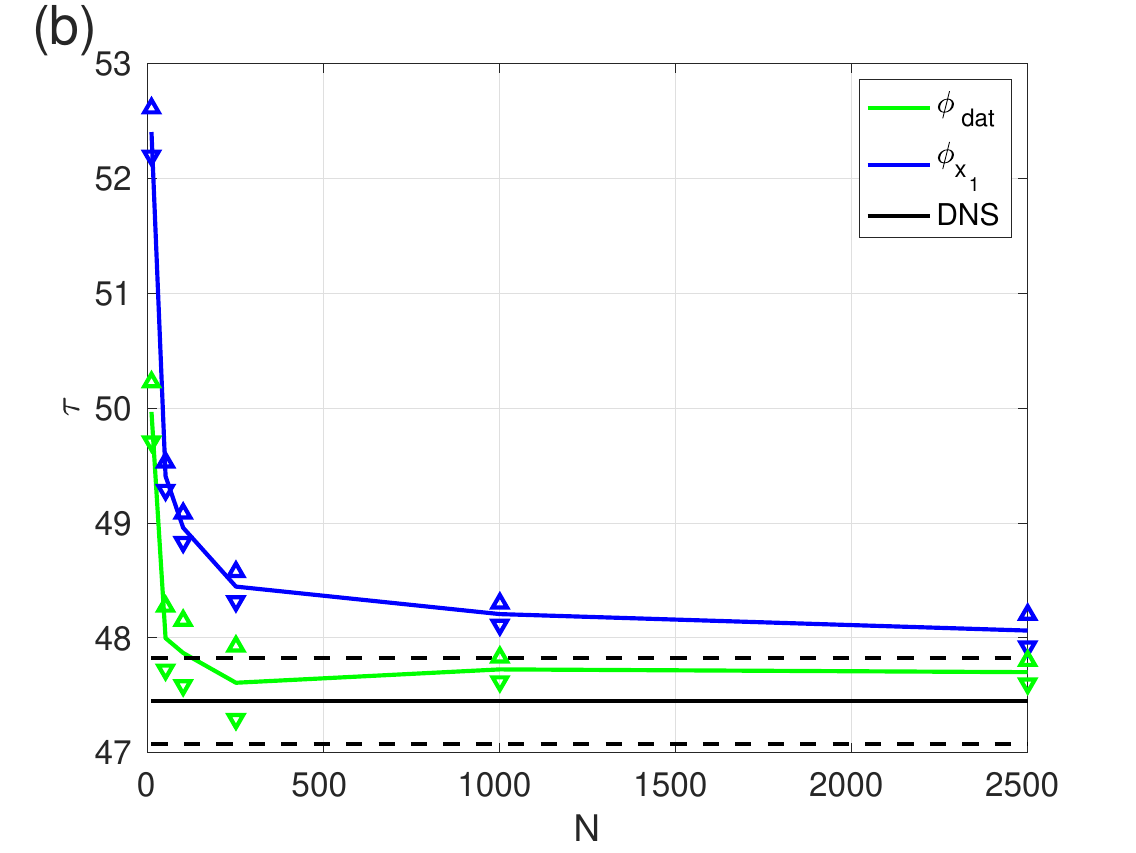}\includegraphics[width=0.33\linewidth]{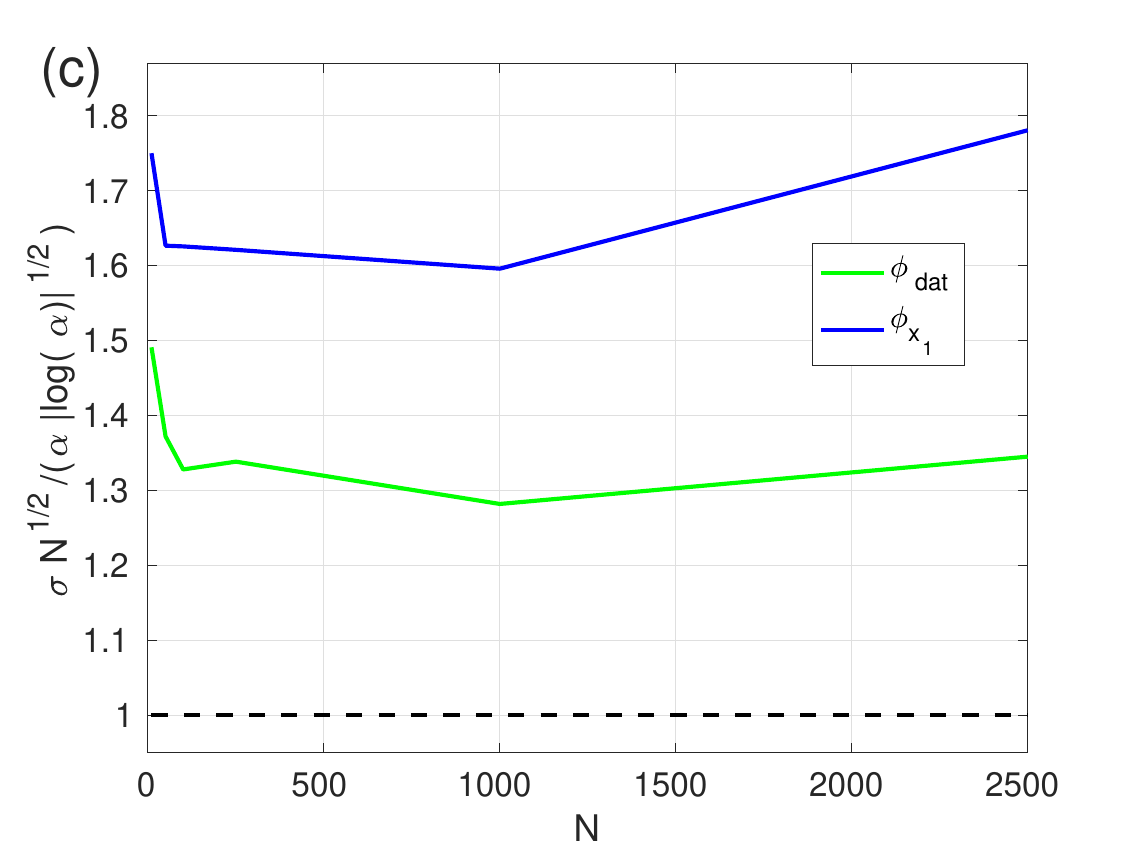}}
\caption{Efficiency of the AMS algorithm with the learned committor function for large dataset as a function of the number of clones $N$, for the Charney--DeVore model.
Comparison of the estimated (a) Transition probability $\langle \alpha\rangle$, (b) Duration of reactive trajectories $\langle \tau\rangle$, and (c) Rescaled standard deviation $\sigma$. For each plot the black curve is the reference one, either the DNS (a and b) or the optimal value 1 (c).
The dashed black lines are the 95\% confidence interval for the DNS. The color curves have been computed using the AMS, with the learned committor function (green) and the linear score function (blue).
The learned committor function gives excellent results, similar to the linear one for the weak apparent bias of the transition probability (a), and much better than the linear one for the standard deviation and the duration of reactive trajectories (b and c).}
\label{fig_cdv}
\end{figure*}
We first show the estimate of the transition probability $\alpha$ as a function of the number of clones used in AMS in figure~\ref{fig_cdv} (a).
For all three estimates, the 95\% intervals of confidence are fairly large: 2\% of the empirical average.
All three intervals overlap if more than $100$ clones are used in AMS computations.
Based on this observable alone, both score functions give comparable results, and we cannot conclude on whether one is better than the other.

We then show the estimate of the average duration of reactive trajectories as a function of the number of clones used in AMS computations in figure~\ref{fig_cdv} (b).
The two AMS estimates $\langle \tau\rangle_M$ decrease with $N$ towards an asymptotic value.
With the learned committor function $\phi_{\rm dat}$, the 95\% confidence intervals of the AMS and DNS estimates overlap if $N\ge 250$.
This never happens for the linear score function $\phi_{x_1}$.

Finally, figure~\ref{fig_cdv} (c) shows the rescaled standard deviation of the AMS estimator of $\alpha$ as a function of the number of clones.
Both are compared to the reference value $1$.
The learned committor function significantly reduces the statistical error, compared to the linear score function.

We conclude that using the learned committor function computed from a long dataset leads to more precise results than using the user-defined score function $\phi_{x_1}$, especially for the statistical error and for the estimate of the duration of reactive trajectories.
We note that AMS computations yield estimates close to the asymptotic value if $N\ge 1000$.

\subsubsection{Efficiency of the AMS algorithm with the learned committor function as a function of the dataset length for the Charney--DeVore model}

As we did with the 2D three-well model (Sec.~\ref{sec:datasize3wells}), we now wish to determine the amount of data necessary to learn a committor function leading to good AMS estimates.
Again, we sample longer and longer trajectories, containing from 1 to 99 transitions.
From each of these datasets we learn a committor function and use it in AMS computations using $N=1000$ clones.
For each dataset length, we perform an average over at least $M=100$ independent realizations of the score function, and up to $M=870$ realisations for the datasets containing 38 transitions.

\begin{figure*}
\centerline{\includegraphics[width=0.33\linewidth]{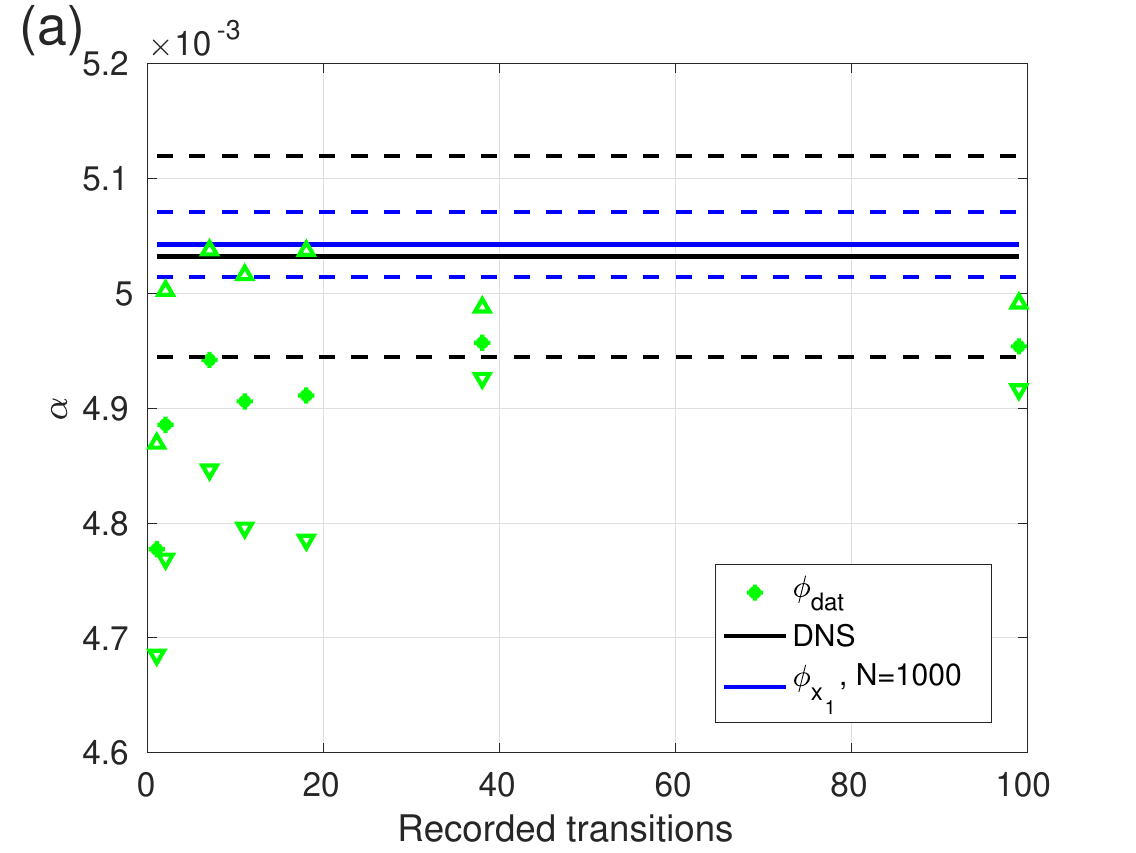}\includegraphics[width=0.33\linewidth]{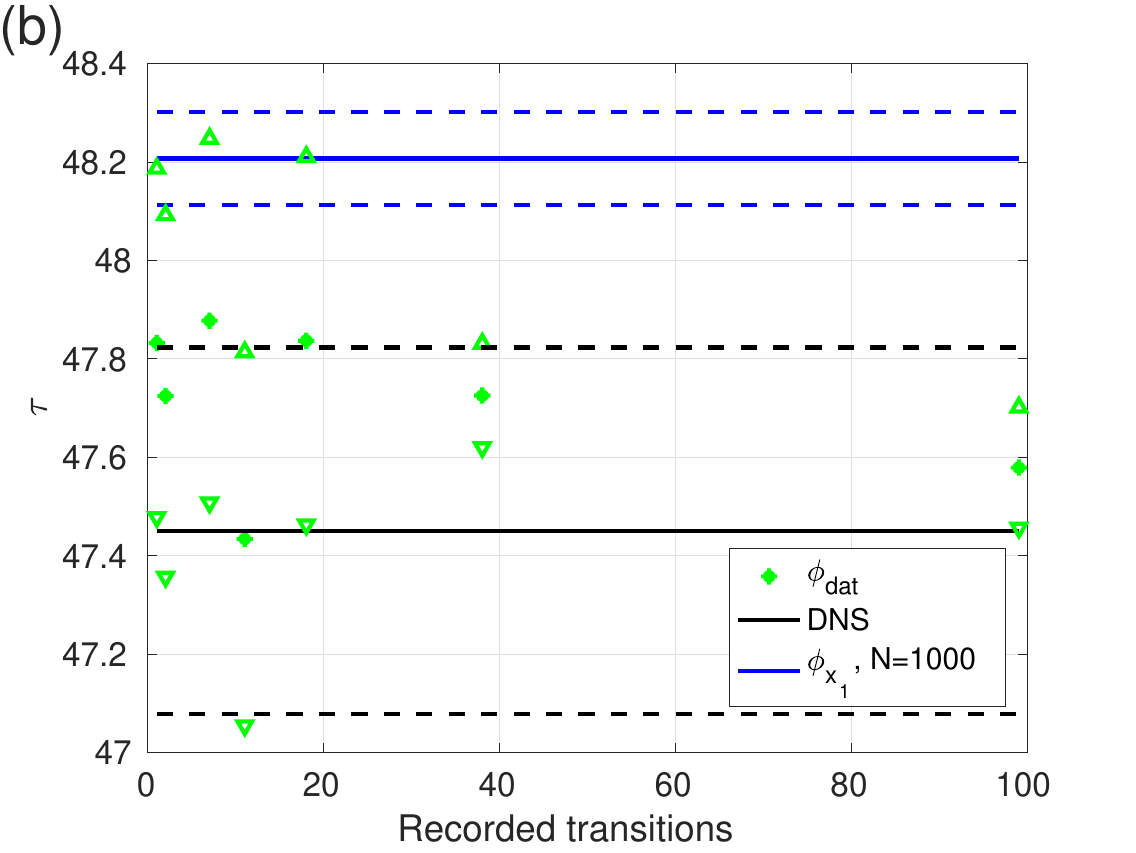}\includegraphics[width=0.33\linewidth]{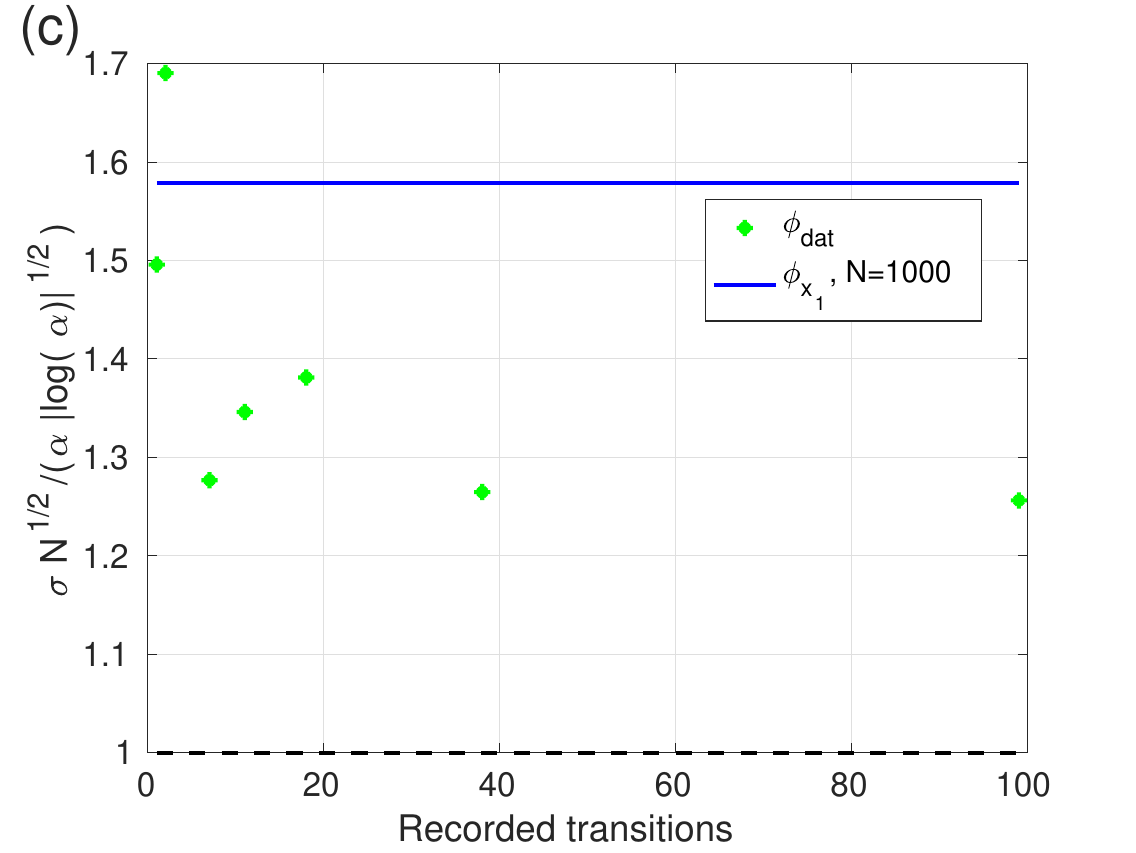}}
\caption{Efficiency of the AMS algorithm with the learned committor function as a function of the dataset length (measured in number of transitions), for the Charney--DeVore model. Comparison of the estimated (a) Transition probability $\langle \alpha\rangle$, (b) Duration of reactive trajectories $\langle \tau\rangle$, and (c) Rescaled standard deviation $\sigma$. For each plot the black curve is the reference one, either the DNS (a and b) or the optimal value 1 (c). The dashed black lines are the 95\% confidence interval for the DNS. The color curves have been computed using the AMS, with the learned committor function (green) and the linear score function (blue). For dataset as short as 5 transitions the AMS algorithm with the learned committor function leads to results as precise as the DNS, and more precise than the linear score function, for both the rescaled standard deviation and trajectory duration.
Having few transitions in the dataset leads to variability in the quality of the score function.}
\label{imp_cdv}
\end{figure*}
We first consider the transition probability $\alpha$ (figure~\ref{imp_cdv} (a)) and the average duration of the reactive trajectories $\tau$ (figure~\ref{imp_cdv} (b)) as a function of the number of recorded transitions.
We note that as soon as there are more than five recorded transitions in the dataset, using the AMS with the learned committor, the 95\% intervals of confidence of $\alpha$ overlap with the DNS estimate, while the confidence interval of $\tau$ always overlap with that of DNS.
This indicates that the learned committor function is relevant for smaller datasets than those used in Sec.~\ref{sec:apparentbiasCDV}.
As the size of the dataset is increased, the confidence interval of DNS estimates and AMS estimates using the learned score function overlap more and more.

We now examine the rescaled variance of the estimate of $\alpha$ performed by AMS, using the learned committor score function, as a function of the number of transitions recorded in the dataset used (figure~\ref{imp_cdv} (c)).
It is compared to the rescaled variance of the estimate of $\alpha$ performed by AMS using the linear score function.
We note that if there are very few transitions recorded in the dataset, the variance is larger than 1.5 and can be larger than the one obtained using the linear score function.
However, this quickly improves with the size of the dataset, the rescaled variance can go down to 1.2 when using the dataset where 38 transitions are recorded.

\section{Conclusions}


In this paper, we have proposed a data-driven approach for the computation of the committor function.
This approach relies on the analogue method to define an effective dynamics starting only from observations.
We have shown that this defines a Markov chain on the observed states of the dynamics, which approximates the true propagator.
A spectral characterization of the committor function can thus be used.
This method of computation of the committor function gives remarkably smooth and robust results for the committor function.

We have highlighted by means of two examples that it is possible to obtain fairly precise estimates of the committor function, even in cases where few observations are available.
In addition, we have pointed out that these approximations are more precise than those provided by a more naive data-driven approach and that increasing the amount of data results in a faster reduction of the error.
The estimates are more precise because the analogue Markov chain is a dynamical approach, which uses all the information contained in the trajectories, while this is not the case for the direct approach, which treats all the points of the same reactive trajectory equally.
We also stress that the analogue Markov chain approach can be used with trajectories of any length, not necessarily distributed according to the invariant measure of the dynamics.

Finally, we provided evidence of the advantage of coupling the analogue method with a rare event algorithm.
Indeed, the learned committor with the analogue Markov chain can be used as a score function performing better than user-defined score functions.
This means that it is possible to develop an almost-fully automatic algorithm that requires very little knowledge and understanding of the system under consideration. The quality of the results suggest that better understanding can be obtained \emph{a posteriori}.
While this work presented one iteration of the loop, we are currently implementing and testing the feedback control loop (figure~\ref{fig:Feedback-Control-Loop}) between the AMS algorithm and the learning of the committor with the analogue Markov chain.

Although the learned committor function based on the analogue Markov chain, and its coupling with rare event algorithms,
have revealed several very interesting advantages, some limitations might arise especially when one faces high-dimensional systems.
We have tested the approach for a fairly complex dynamics with 6 degrees of freedom. It still has to be tested for more complex dynamics.
For systems in high dimensions, the choice of the distance for the analogue method might be a critical issue.

Another interesting question would be to compare the quality of the estimation of the committor function using the analogue Markov chain with other methods.
It could be compared to other methods based on dynamical information, sometimes more complicated, for instance the direct Galerkin approximation~\cite{thiede2019galerkin,strahan2021long} method.
It would also be interesting to compare it to direct approaches using machine learning.

\ack
The work of D.~Lucente was funded through the ACADEMICS grant of IDEXLYON, project of the Université de Lyon, PIA operated by ANR-16-IDEX-0005. This work was supported by the ANR grant SAMPRACE, project ANR-20-CE01-0008-01. This publication was supported by a Subagreement from the Johns Hopkins University with funds provided by Grant No.~663054 from Simons Foundation. Its contents are solely the responsibility of the authors and do not necessarily represent the official views of Simons Foundation or the Johns Hopkins University.  Computer time was provided by the \enquote{P\^ole Scientifique de Mod\'elisation Num\'erique} and the \enquote{Centre Blaise Pascal}  in Lyon, as well as by the computing resources made available by the LMFL in Lille.
The authors thank Charles-\'Edouard Br\'ehier, Arnaud Guyader, Tony Leli\`evre, Cl\'ement Le Priol, Eric Simonnet and Pascal Yiou for helpful comments on the manuscript.

\section*{References}

\bibliographystyle{unsrt}
\bibliography{references_analogue.bib}

\end{document}